\documentclass[a4paper, 11pt]{article} 
\usepackage{amsfonts,amsmath,amssymb,amscd,mathabx}
\usepackage{latexsym}
\usepackage{graphicx}
\usepackage{color} 
\usepackage{hyperref}

\newcommand{\dxp}{d_{|x} p}
\newcommand{\spd}{\sqrt{\pi \delta_L  d^{n+1}}}
\newcommand{\nax}{\nabla_{|x}}
\newcommand{\elstar}{(E\otimes L^d)^*}
\newcommand{\kxp}{\ker d_{|x} p}
\newcommand{\pix}{\pi_1 (\pi_2^{-1}(x))}
\newcommand{\invpin}{\frac{1}{\sqrt \pi^n } }

\newcommand{\rdp}{\R \Delta^d_p}
\newcommand{\vkn}{V_{k-1,n-1}}

\newcommand{\mkr}{M_{k-1}(\R)}
\newcommand{\lhp}{L(H,H^\perp)}
\newcommand{\grkn}{\text{Gr}(k-1,n-1)}
\newcommand{\trace}{\text{Tr }}

\newcommand{\bh}{\mathcal H}

\newcommand{\Nn}{\mathbb{N}}
\newcommand{\R}{\mathbb{R}}
\newcommand{\E}{\mathbb{E}}

\newcommand{\Zz}{\mathbb{Z}}
\newcommand{\C}{\mathbb{C}}

\newcommand{\ep}{\epsilon}
\newcommand{\si}{\sigma}
\newcommand{\Si}{\Sigma}

\newcommand{\la}{\lambda}

\newcommand{\dbar}{\bar \partial}
\newcommand{\up}{{(U,P)}}
\newcommand{\ups}{{(U_S,P_S)}}

\newcommand{\beqr}{\begin{eqnarray}}
\newcommand{\eeqr}{\end{eqnarray}}

\newcommand{\beq}{\begin{eqnarray*}}
\newcommand{\eeq}{\end{eqnarray*}}

\newcommand{\bq}{\begin{equation}}
\newcommand{\eq}{\end{equation}}

\newcommand{\bpr}{\begin{preuve}}
\newcommand{\epr}{\end{preuve}}

\newenvironment{preuve}[1][]
{\vskip 2mm  {\it \bf Proof#1. }}{$\Box$ \vskip 2mm}

\newtheorem{Theorem}{Theorem}[subsection]
\newtheorem{Remark}[Theorem]{Remark}
\newtheorem{Definition}[Theorem]{Definition}
\newtheorem{Lemma}[Theorem]{Lemma}
\newtheorem{Proposition}[Theorem]{Proposition}
\newtheorem{Corollary}[Theorem]{Corollary}

\newcommand{\equid}{\underset{d\to \infty}{\sim}}

\newcommand{\eld}{E\otimes L^d}
\newcommand{\hxed}{H^0(X,E\otimes L^d)}
\newcommand{\rhxed}{\R H^0(X,E\otimes L^d)}
\newcommand{\held}{h_{E,d }}

\newcommand{\csi}{C_\si}
\newcommand{\rcsi}{\R C_\si}

\newcommand{\prcsi}{p_{|\R C_\sigma}  }

\newcommand{\cpun}{\C P^1}

\newcommand{\lap}{ \dot{\lambda}}
\newcommand{\lapo}{ \dot{\lambda}_0}
\newcommand{\sip}{ \dot{\si}}
\newcommand{\vp}{ \dot{v}}
\newcommand{\vpo}{\dot{v}_0}

\newcommand{\xp}{\dot{x}}

\title{Expected topology of random \\ real algebraic submanifolds}
\author{Damien Gayet, Jean-Yves Welschinger}

\begin{document}
\large
\maketitle
\centerline{\textbf{Abstract}}
Let $X$ be a smooth complex projective manifold of dimension $n$ equipped with an ample line bundle $L$ and a rank $k$ holomorphic vector bundle $E$. 
We assume that $1\leq k \leq n$,  that $X$, $E$ and $L$ are defined over the reals
and denote by $\R X$ the real locus of $X$. Then, we estimate from above and below
the  expected Betti numbers of the vanishing loci in $\R X$ 
 of holomorphic real sections of
$E\otimes L^d$,  where $d$ is a large enough integer. Moreover,
given any closed connected codimension $k$ submanifold 
$\Sigma$ of $\R^n$ with trivial normal bundle,
 we prove that a real section of $E\otimes L^d$ has a positive probability, independent of $d$,
 to contain around $\sqrt d^n $ connected components diffeomorphic to $\Sigma$ in its vanishing locus. \\


\tableofcontents

\section{Introduction}
Let $X$ be a smooth complex projective manifold of positive dimension $n$ equipped with an ample line bundle $L$ and let $E$ be a holomorphic vector bundle of rank $k$ over $X$. 
From the vanishing theorem of Kodaira and Serre, 
we know that 
the dimension $N_d$ of the complex vector space $H^0(X,\eld)$ of global holomorphic sections of $\eld$ 
%
grows
as a polynomial of degree $n$ in $d$. 
We will assume throughout this paper that $1\leq k \leq n$ and that $X$, $E$ and $L$ are defined over the reals. 
We denote by $\R X$ the real locus of $X$ and by $\rhxed$ the real vector space  of real holomorphic sections of $\eld$, see (\ref{rsections}). Its dimension equals $N_d$. 
The discriminant locus $\R \Delta_d\subset \rhxed$ of sections which do not vanish transversally 
is a codimension one submanifold for $d$ large enough and for every $\si$ in its complement, the real vanishing locus $\rcsi $ of $\si$ is a smooth
codimension $k$ submanifold of $\R X$. The topology of $\rcsi$ drastically depends on the choice of $\si \in \rhxed\setminus \R \Delta_d$. 
When $n=k=1$, $X=\C P^1$, $L = \mathcal O_{\C P^1} (1)$ and $E=\mathcal O_{\C P^1}$ for example, $\si$ is  a real polynomial of degree $d$ in one variable and $\rcsi$ the set of its real roots. 

The space $\rhxed$ inherits classical probability measures. Indeed, let $h_E$ be a Hermitian metric on $E$ and $h_L$ be a Hermitian metric of positive curvature on $L$, both $h_E$ and $h_L$ being real, that is invariant under the $\Zz/2\Zz$-Galois action of $E$ and $L$. We denote by $\held = h_E\otimes h_L^d$ the induced metric on $\eld$. Then, 
the vector space $\rhxed$ becomes Euclidean, with the $L^2$-scalar product defined by 
$$ \forall \si, \tau \in \rhxed, \ \langle \si, \tau \rangle = 
\int_X \held (\si, \tau) dx,$$
where $dx$ denotes any chosen volume form on $X$ (our results being asymptotic in $d$, they 
 turn out not to depend on the choice of $dx$). 
It thus inherits a Gaussian probability measure $\mu_\R$ 
whose density at $\si \in \rhxed$ with respect to the Lebesgue measure is
$\frac{1}{\sqrt \pi^{N_d}} e^{-\|\si\|^2}.$

What is the typical topology of $\rcsi$ for $\si \in \rhxed$
 chosen at random for $d\mu_\R$? We do not know, but can estimate its average Betti numbers. 
 To formulate our results, let us denote, for every $i\in \{0, \cdots, n-k\}$, 
 by $b_i(\rcsi, \R )= \dim H_i(\rcsi, \R)$ the $i$-th Betti number of $\rcsi$ 
and by 
$$ \E(b_i) = \int_{\rhxed \setminus \R \Delta_d} b_i(\rcsi, \R) d\mu_\R (\si)$$ 
its expected value. 
\subsection{Upper estimates}
As in \cite{GaWe3}, for every $i\in \{0, \cdots, n-k\}$, we denote by $Sym_\R (i,n-k-i)$ the open
cone of real symmetric matrices of size $n-k$ and signature $(i,n-k-i)$, by $\mu_\R$ the 
classical Gaussian measure on the space of real symmetric matrices and by 
$e_\R (i,n-k-i)$ the numbers 
\bq\label{ern}
e_\R (i,n-k-i) = \int_{Sym_\R (i,n-k-i)} |\det A| d\mu_\R (A),
\eq
see \S \ref{statement}.
We then denote by $Vol_{h_L} (\R X) $ the volume of $\R X$ for the Riemannian metric induced
by the K\"ahler metric $g_{h_L}$ defined by the curvature form of  $h_L$, see (\ref{omega}) and (\ref{metrique}).
\begin{Theorem}\label{th11}
Let $X$ be a smooth real projective manifold of dimension $n$, $(L,h_L)$ be
a real holomorphic Hermitian line bundle of positive curvature over $X$ and $(E,h_E)$
be a rank $k$ real holomorphic Hermitian vector bundle, with $1\leq k \leq n$, $k\neq n$. Then, for every $0\leq i\leq n-k$,
\beq
\limsup_{d\to \infty} \frac{1}{\sqrt d^n } \E(b_i)
& \leq & {n-1 \choose k-1} e_\R (i,n-k-i) \frac{Vol_{h_L} (\R  X)}{Vol_{FS}(\R P^k)}.
\eeq
Moreover, when $k=n$, $\frac{1}{\sqrt d^n } \E(b_0) $
converges to $\frac{Vol_{h_L} (\R  X)}{Vol_{FS}(\R P^n)}$ 
as $d$ grows to infinity.
\end{Theorem}
In fact, the right hand side of the inequality given by Theorem \ref{th11} 
also involves the determinant of random matrices of size $k-1$
and the volume of the Grassmann manifold of $(k-1)$ linear subspaces of $\R^{n-1}$, 
see Theorem \ref{th1}, but these can be computed explicitly. 
Note that when $E$ is the trivial line bundle, Theorem \ref{th11} reduces to Theorem 1.1 of \cite{GaWe3}.

Theorem \ref{th11} relies on Theorem \ref{theo 3}, which 
establishes the asymptotic equidistribution of clouds of critical points, see \S \ref{statement}. 
We obtain a similar result in a complex projective setting, for critical points of Lefschetz pencils,
see Theorem \ref{lefschetz}.

\subsection{Lower estimates and topology}\label{intro lower}
Let $\Sigma$ be a closed submanifold of codimension $k$  of $\R^n$, $1\leq k \leq n$,
which we do not assume to be connected. 
For every $\si \in \rhxed\setminus \R \Delta_d$,
we denote by 
$N_\Sigma(\si) $ the maximal number of disjoint open subsets 
of $\R X$ having the property that each such open subset $U'$ 
contains a codimension $k$ submanifold $\Sigma'$
such that $\Sigma'\subset \rcsi$
and $(U',\Sigma')$
is diffeomorphic to $(\R^n , \Sigma)$. We then set 
\bq\label{nsigma}
 \E(N_\Sigma)= \int_{ \rhxed\setminus \R \Delta_d} N_\Sigma (\si) d\mu_\R(\si)
 \eq
and we associate to $\Sigma$, in fact to its isotopy class in $\R^n$,
a constant $c_\Sigma$ which is positive if and only if $\Sigma$ has trivial normal bundle in $\R^n$, see (\ref{csigma})
for its definition and Lemma \ref{lemmaN}. The latter measures \`a la Donaldson the amount of transversality that a polynomial
map $\R^n \to \R^k$ vanishing along a submanifold isotopic to $\Sigma$ may have. 
\begin{Theorem}\label{th2}
Let $X$ be a smooth real projective manifold of dimension $n$, $(L,h_L)$ be
a real holomorphic Hermitian line bundle of positive curvature over $X$ and $(E,h_E)$
be a rank $k$ real holomorphic Hermitian vector bundle, with $1\leq k \leq n$.  Let $\Sigma$
be a closed submanifold of codimension $k$ of $\R^n$ with trivial normal bundle, which does not 
need to be connected. Then, $$
\liminf_{d\to \infty} \frac{1}{\sqrt d^n} \E(N_\Sigma)\geq c_\Sigma Vol_{h_L}(\R X).$$
\end{Theorem}

In particular, when $\Sigma$ is connected, Theorem \ref{th2} 
bounds from below the expected number of connected 
components diffeomorphic to $\Sigma$
in the real vanishing locus 
of a random section $\si\in \rhxed$. The constant $c_\Sigma$ 
does not depend on the choice of
the triple $(X, (L,h_L), (E,h_E)), $ it only depends on $\Sigma$.
When $k = 1$ and $E = \mathcal O_X$, 
Theorem \ref{th2} coincides with Theorem 1.2 of \cite{GaWe4}. 
Computing $c_\Sigma$ for explicit submanifolds $\Sigma$ yields  the following lower bounds for the Betti numbers.
\begin{Corollary}\label{coro i}
Under the hypotheses of Theorem \ref{th2}, 
for every $i\in \{0, \cdots, n-k\}$, 
$$ \liminf_{d\to \infty} \frac{1}{\sqrt d^n} \E(b_i )\geq \exp(-e^{84 +6n}) Vol_{h_L}(\R X).$$
\end{Corollary}

\subsection{Some related results}
The case $X= \C P^1$, $E=\mathcal O_{\C P^1}$ and $L = \mathcal O_{\C P^1}(1)$
 was first considered by M. Kac in  \cite{Kac} for a different measure. 
In this case and with our measure,  Kostlan \cite{Ko} and Shub and Smale  \cite{SS} gave an exact formula for the mean number of real roots 
of  a polynomial, as well as  the mean number of intersection points of $n$  hypersurfaces in $\R P^n$. Still in $\R P^n$, Podkorytov 
\cite{Podkorytov} computed the mean Euler characteristics of random algebraic hypersufaces, and B\"urgisser \cite{Burgisser} 
extented this result to complete intersections. In \cite{GaWe1},
we proved the exponential rarefaction of real curves 
with a maximal number of components in real algebraic surfaces. In \cite{GaWe2} and \cite{GaWe3},
we bounded from above the mean Betti numbers of random real hypersurfaces in
real projective manifolds and in \cite{GaWe4}, we gave a lower bound  for them. 

A similar probabilistic study of complex projective manifolds has been performed
by Shiffman and Zelditch, see \cite{SZ}, \cite{SZ2}, \cite{BSZ2} for example, 
or also \cite{BBL}, \cite{Sodin-Tsirelson}. In particular,
the asymptotic equidistribution 
of critical points
of random sections over a fixed  projective manifold
has been studied in
\cite{DSZ1}, \cite{DSZ2} and  \cite{MacD}, or also  \cite{Fyodorov}, \cite{ABC}, \cite{DeMa},
while   we 
studied 
 critical points of the restriction of a fixed 
Morse function on random real hypersurfaces, see \cite{GaWe2}, \cite{GaWe3}.

A similar question concerns
the mean number of components of the vanishing locus of  eigenfunctions of the Laplacian.
It has been studied on the round sphere by Nazarov and Sodin \cite{Nazarov-Sodin} (see also \cite{Sodin}),
Lerario and Lundberg \cite{Lerario-Lundberg} or Sarnak and Wigman \cite{SarnakWigman}.  
In a general Riemannian setting, Zelditch proved in \cite{Zelditch} the
equidistribution of the vanishing locus,  
  whereas critical points of random eigenfunctions of the Laplacian 
have been addressed by Nicolaescu in \cite{Nicolaescu}. \\

Section \ref{partie 2} is devoted to lower estimates and the proof of Theorem \ref{th2}. In this proof, the $L^2-$estimates
of H\"ormander play a crucial r\^ole, see \S \ref{para II.1},
and we follow the same approach as in \cite{GaWe4} (see also \cite{Gayet} for a similar construction).
Section \ref{II} is devoted to upper estimates and the proof of Theorem \ref{th11}. \\

\textit{Aknowledgements.} The research leading to these results has received funding
from the European Community's Seventh Framework Progamme 
([FP7/2007-2013] [FP7/2007-2011]) under
grant agreement $\text{n}\textsuperscript{o}$ [258204].

\section{Lower estimates for the expected Betti numbers}\label{partie 2}

\subsection{Statement of the results} \label{I.1}
\subsubsection{Framework}

Let us first recall our framework. We denote by $X$
a smooth complex projective manifold of dimension $n$
defined over the reals, by $c_X: X \to X$ the induced Galois antiholomorphic involution 
 and by $\R X= \text{Fix} (c_X)$ 
the real locus of $X$ which we implicitly  assume to be non-empty.
We then consider an ample line bundle $L$ over $X$,
also defined over the reals. It comes thus equipped with an antiholomorphic 
involution $c_L: L \to L$ which turns the bundle projection map $\pi : L\to X$
into a $\Zz /2\Zz $-equivariant one, so that $c_X \circ \pi = \pi \circ c_L$.
We equip $L$ in addition with a real Hermitian metric $h_L$, thus 
invariant under $c_L$,
 which has a positive curvature form $\omega$ locally defined by 
 \bq \label{omega}
 \omega = \frac{1}{2i\pi}
\partial \dbar \log h_L(e,e) 
\eq
 for any  non-vanishing local
holomorphic section $e$ of $L$. This metric induces a K\"ahler 
metric 
\bq \label{metrique}
g_{h_L}= \omega(.\, ,i.\, )
\eq
 on $X$, which reduces to
a Riemannian metric $g_{h_L}$ on $\R X$.
Let finally $E$ be a holomorphic vector bundle of rank $k$, $1\leq k\leq n$,
defined over the reals and equipped with a antiholomorphic involution $c_E$ and a real Hermitian metric $h_E$. 
For every $d>0$, we denote by 
\bq \label{rsections}
\rhxed = \{\si \in \hxed\ | \ (c_{E}\otimes c_{L^d}) \circ \si = \si \circ c_X\}
\eq
the space of global real holomorphic sections of $E\otimes L^d$. It is equipped with the $L^2$-scalar
product defined by the formula
\begin{equation}\label{l2product}
 \forall (\si, \tau)\in \rhxed, \ \langle \si, \tau\rangle = \int_X h_{E,d} (\si,\tau)(x) dx, 
\end{equation}
where $h_{E,d}= h_E\otimes h^d_L$. Here,  $dx$
denotes any  volume form of $X$. 
For instance, $dx$ can be chosen to be the normalized volume form $dV_{h_L}= \frac{\omega^n }{\int_X \omega^n }$.
This $L^2$-scalar 
product finally induces a Gaussian probability measure $\mu_\R$
on $\rhxed$
whose density with respect to the Lebesgue one at $\si\in \rhxed$ 
writes $\frac{1}{\sqrt\pi^{N_d}}e^{-\|\si\|^2}$,
where $N_d= \dim \hxed.$
It is with respect to this probability measure that we consider random real codimension $k$ submanifolds (as
in the works \cite{Ko} and \cite{SS}, \cite{GaWe2}, \cite{GaWe3} and \cite{GaWe4}). 

\subsubsection{The lower estimates}
The aim of Section \ref{partie 2} is to prove Theorem \ref{th2}.
In addition to Theorem \ref{th2}, we also get the following 
Theorem \ref{theo 0}, which is a consequence of Proposition \ref{prop 5} below.
\begin{Theorem}\label{theo 0}
Under the hypotheses of Theorem \ref{th2}, for every $0\leq \epsilon <1$, 
$$ \liminf_{d\to \infty } \mu_\R \left\{ \si \in \rhxed \ | \ N_\Sigma (\sigma )\geq \epsilon c_\Sigma Vol_{h_L}(\R X) \sqrt d^n \right\} >0.$$
\end{Theorem}
In fact, the positive lower bound given by Theorem \ref{theo 0} can be made  explicit, see
(\ref{minoration}). 

Let us now denote, for every  $1\leq k\leq n$, by $\bh_{n,k}$ the set of diffeormophism
classes of smooth closed  connected codimension $k$ submanifolds of  $\R^n$. For every $i\in \{ 0, \cdots, n-k\}$
and every $[\Sigma]\in \bh_{n,k}$, we denote by $b_i(\Sigma) =\dim H_i(\Sigma ;\R)$
its $i$-th Betti number with real coefficients and
by $m_i(\Sigma)$
its $i$-th Morse number. This is the infimum over all Morse functions $f$ on $\Sigma$
of the number of critical points of index $i$ of $f$. 
 Then, 
we set $c_{[\Sigma]}= \sup_{\Sigma \in [\Sigma]  } c_\Sigma $ and
$$ \E(m_i) = \int_{\rhxed \setminus \R \Delta_d} m_i(\rcsi) d\mu_\R(\si) .$$
\begin{Corollary}\label{coro 1}
Let $X$ be a smooth real projective manifold of dimension $n$, $(L,h_L)$ be
a real holomorphic Hermitian line bundle of positive curvature over $X$ and $(E,h_E)$
be a rank $k$ real holomorphic Hermitian vector bundle, with $1\leq k \leq n$. 
Then,
for every $i\in\{0, \cdots, n-k\}$, 
\begin{eqnarray}
\liminf_{d\to \infty} \frac{1}{\sqrt d^n} \E(b_i) &\geq &\big(\sum_{[\Sigma]\in \bh_{n,k}} c_{[\Sigma]} b_i(\Sigma) \big)Vol_{h_L}(\R X)  \label{loest}
\text{ and likewise }\\
\liminf_{d\to \infty} \frac{1}{\sqrt d^n} \E(m_i)&\geq &\big(\sum_{[\Sigma]\in \bh_{n,k}} c_{[\Sigma]} m_i(\Sigma)\big) Vol_{h_L}(\R X). \label{loest'}
\end{eqnarray}

\end{Corollary}

Note that in Corollary \ref{coro 1}, we could have chosen one representative $\Sigma$ in each diffeomorphism class $[\Sigma]\in \bh_{n,k}$
and obtained the lower estimates (\ref{loest}), (\ref{loest'}) with constants $c_\Sigma$ instead of $c_{[\Sigma]}$. But it turns out that in the proof of Corollary \ref{coro 1}
we are free to choose the representative we wish in every diffeomorphism class and that the higher $c_\Sigma$ is, the better the estimates (\ref{loest}), (\ref{loest'}) are. 
This is why we introduce the constant $c_{[\Sigma]}$, which is positive if and only if $[\Sigma]$ has a representative $\Sigma$ with trivial normal bundle in
$\R^n$, see   (\ref{csigma}) and Lemma \ref{lemmaN}.

\subsection{Closed affine real algebraic submanifolds}\label{para I.2}

We introduce here the notion of regular pair, see Definition \ref{def 1},
and the constant $c_\Sigma$ associated to any isotopy class of smooth closed codimension $k$ submanifold $\Sigma$
of $\R^n$, see (\ref{csigma}).

\begin{Definition}\label{def 1}
Let $U$ be a bounded open subset of $\R^n$ and $P\in \R[x_1, \cdots x_n]^k$, $1\leq k \leq n$. 
The pair $(U,P)$ is said to be regular if and only if 
\begin{enumerate}
\item zero is a regular value of the restriction of $P$ to  $U$,
\item the vanishing locus of $P$ in $U$ is compact.
\end{enumerate}
\end{Definition}

Hence, for every regular pair $(U,P)$, the vanishing locus of $P$ does not intersect the boundary of $U$ and it meets $U$
in a smooth compact codimension $k$ submanifold. 

In the sequel, for every integer $p$ and every vector $v\in \R^p$, we denote by $|v|$ its Euclidian norm, and for every 
 integers $p$ and $q$, and every  linear map $F : \R^p \to \R^q$,
we denote by $F^*$ the adjoint of $F$, defined by the property
$$\forall v\in \R^p, \forall w\in \R^q,\ 
\langle F(v), w\rangle  = \langle v,  F^*(w)\rangle,$$
and denote by  $\| F\|$ its operator norm, that is 
$$ 
 \| F \| = \sup_{v\in \R^p\setminus \{0\}} | F(v) | /|v|.
 $$
 We will also use the norm 
 $$\|F\|_2 = \sqrt{\trace FF^*}.$$
 These norms satisfy $\|F\| \leq \|F\|_2.$
Finally, if $P = (P_1, \cdots, P_k) \in \R [x_1, \cdots, x_n]^k$, we 
 denote by $ \|P\|_{L^2}$ its $L^2$-norm defined by 
 \begin{equation}\label{polyL2}
  \|P\|^2_{L^2} = \int_{\C^n }|P(z)|^2 e^{-\pi|z|^2} dz = \sum_{i=1}^k \int_{\C^n }|P_i(z)|^2 e^{-\pi |z|^2} dz = \sum_{i=1}^k \|P_i\|^2_{L^2}.
  \end{equation}

\begin{Definition}\label{def 2}
For every regular pair $(U,P) $ given by Definition \ref{def 1}, 
we denote by $\mathcal T_{(U,P)} $ the set of $(\delta, \epsilon)\in (\R^*_+)^2$ such that
\begin{enumerate}
\item there exists a compact subset $K$ of $U$ satisfying $\inf_{x\in U\setminus K} |P(x)| >\delta$,
\item for every $y\in U$, $|P(y)| <\delta \Rightarrow 
\forall w\in \R^k, \ |(d_{|y}P)^*(w) |\geq \epsilon |w|$.
\end{enumerate}
\end{Definition}
Hence, for every regular pair $(U,P)$ given by Definition \ref{def 1}, $(\delta, \epsilon)$ belongs to $\mathcal T_{(U,P)} $ provided the 
$\delta$-sublevel of $P$ does not intersect the boundary of $U$ while inside this $\delta$-sublevel, $P$ is in a sense $\epsilon$-far
from having a critical point. This quantifies how much transversally $P$ vanishes in a way similar to the one used by Donaldson in \cite{Donaldson}.

Then, for every regular pair $\up$, we set $R_{(U,P)} = \max (1, \sup_{y\in U} |y|) $, so that $U$ is contained in the ball centered at the origin and of radius
$R_{(U,P)}$. Finally, we set
\bq \label{tauUP}
\tau_{(U,P)} = 24 k\rho_{R_{(U,P)}} \|P\|^2_{L^2} \inf_{(\delta,\epsilon)\in {\mathcal T_{(U,P)}}} (\frac{1}{\delta^2}+ \frac{\pi n}{\epsilon^2})\in \R^*_+,
\eq
where, for every $R>0$, 
\bq\label{rhoR}
\rho_R = \inf_{\R^+} g_R,
\eq
\bq\label{gr}
g_R :
s \in \R^*_+ \mapsto \frac{(R+s)^{2n}}{s^{2n}} e^{\pi (R+s)^2},
\eq
so that 
\bq\label{estimation r}
e^{\pi R^2 }\leq \rho_R \leq 4^n e^{4\pi R^2}.
\eq
This constant $\tau_{(U,P)}$ is the main ingredient in the definition of $c_\Sigma $, see (\ref{csigma}). The lower $\tau_{(U,P)}$ is, the larger
$c_\Sigma $ is and the better the estimates given by Theorem \ref{th2} are. Note that $\tau_{(U,P)}$ remains small whenever $\delta, \epsilon$ are not too small,
that is when $P$ vanishes quite transversally in $U$.

Now, let $\Sigma$ be a closed submanifold of codimension $k$  of $\R^n$, not necessarily connected. 
We denote by $\mathcal I _\Sigma$ the set of regular pairs $(U,P)$ given by Definition \ref{def 1}, such that
the vanishing locus of $P$ in $U$ contains a subset isotopic to $\Sigma$ in $\R^n$.

\begin{Lemma}
\label{lemmaN}
Let $\Sigma$ be a closed submanifold of codimension $k >0$  of $\R^n$, not necessarily connected. 
Then, $\mathcal I_\Sigma$ is non empty if and only if the normal bundle of $\Sigma$ in $\R^n$ is trivial. 
\end{Lemma}
\bpr
If $(U,P) \in \mathcal I _\Sigma$, then $P : \R^n \to \R^k$ contains in its vanishing locus a codimension $k$ submanifold
$\widehat{\Sigma}$ which is isotopic to $\Sigma$ in $\R^n$. The normal bundle of $\Sigma$ in $\R^n$ is thus trivial if and 
only if the normal bundle of $\widehat{\Sigma}$ in $\R^n$ is trivial. But the differential of $P$ at every point of $\widehat{\Sigma}$ 
provides an isomorphism between the normal bundle of $\widehat{\Sigma}$ in $\R^n$ and the product $\widehat{\Sigma} \times \R^k$.

Conversely, if $\Sigma$ has a trivial normal bundle in $\R^n$, it has been proved by Seifert \cite{Seif} (see also \cite{Nash}) that there exist
a polynomial map $P : \R^n \to \R^k$  and a tubular neighbourhood $U$ of $\Sigma$ in $\R^n$ such that $P^{-1} (0) \cap U$ is isotopic to 
$\Sigma$ in $U$. The strategy of the proof is to first find a smooth function $U \to \R^k$ in a neighborhood of $\Sigma$ which vanishes transversally
along $\Sigma$ and then to suitably approximate the coordinates of this function by some polynomial, see \cite{Seif}, \cite{Nash}.
The pair $(U,P)$ then belongs to $\mathcal I_\Sigma$ by Definition \ref{def 1}.
\epr

We then set $c_\Sigma = 0$ if $\Sigma$ does not have a trivial normal bundle in $\R^n$ and 
\bq\label{csigma}
c_\Sigma = \sup_{(U,P)\in \mathcal I_\Sigma} \left(\frac{m_{\tau_{(U,P)}}}{2^n Vol(B(R_{(U,P)}))}\right) \text{ otherwise},
\eq
where  $Vol(B(R_{(U,P)}))$ denotes the  volume of the Euclidean ball of radius $R_{(U,P)}$ in $\R^n$,
and where, for every $\tau>0$,
\bq\label{mtau}
 m_\tau = \sup_{ [\sqrt \tau , +\infty[}f_\tau,
\eq
with 
$
 f_\tau :
 a \in  [\sqrt \tau , +\infty[ \  \mapsto \frac{1}{\sqrt \pi}(1-\frac{\tau}{a^2}) \int_a^{+\infty }e^{-t^2} dt.
 $ 
For large values of $m_\tau$, as the ones which appear in \S \ref{para exemples}, the estimate
\bq \label{estimation c} 
c_\Sigma \geq e^{-2 \tau_{(U,P)}}
\eq
holds, compare ($2.8$) of \cite{GaWe4}.

\subsection{H\"ormander sections}\label{para II.1}
Our key tool to prove Theorems \ref{th11} and \ref{th2}
has been developped by L. H\"ormander. We introduce in this \S \ref{para II.1}
the material we need. 
For every positive $d$ and every $\si \in \rhxed$, we set 
$$ \| \si \|^2_{L^2 (h_L)} = \int_X \|\si\|^2 _{h_{E,d}} dV_{h_L},$$
where $dV_{h_L}= {\omega^n}/\int_X \omega^n $, compare (\ref{l2product}).
Let us choose a field of $h_L$-trivializations of $L$ on $\R X$ given by Definition
4 of \cite{GaWe4}. It provides in particular, for every $x\in \R X$,
a local holomorphic chart $\psi_x : (W_x,x) \subset X \to (V_x,0)\subset \C^n$
isometric at $x$, and a non-vanishing holomorphic section $e$ of $L$ defined over  $W_x$ such that
$\phi = -\log h_L(e,e)$ vanishes at $x$ and is positive elsewhere. Moreover,
 there exist a positive constant $\alpha_1$ such that 
\begin{equation}\label{phi}
\forall y\in V_x,  \vert  \phi\circ \psi_x^{-1}(y)- \pi |y|^2 \vert \leq \alpha_1 | y|^3.
 \end{equation}
Restricting $W_x$ if necessary, we choose a holomophic trivialization $(e_1, \cdots, e_k)$ of $E_{|W_x}$ which
is orthonormal at $x$. 
This
provides a trivialization $(e_1\otimes e^d, \cdots, e_k \otimes e^d)$ of $E\otimes L^d_{|W_x}$. 
In this trivialization, the restriction of $\si$ to $W_x$
writes 
\bq\label{sif}
\si = \sum_{j=1}^k f^j_\si e_j \otimes e^d
\eq
 for some holomorphic functions
$f^j_\si : W_x \to \C$,
  We write $f_\si= (f^1_\si, \cdots, f^k_\si)$
  and we set
  \bq \label{bal}
  |\si| = |f_\si|, 
  \eq 
   so that on $W_x$,
 $ \| \si \|^2_{\held} = \big\| \sum_{j=1}^k f_\si^k e_j \big\|^2_{h_E} e^{-d\phi}$  and
 $\|\si(x)\|^2_{\held} = |\si(x)|^2$ since the frames $(e_1, \cdots, e_k )$ and $e$
 are orthonormal at the point $x$ so that in particular $\phi (x) = 0$. For every $z\in W_x$, we define
 \bq \label{bal2}
  \|d_{|z}\si\|_2= \|d_{|y} (f_\si \circ \psi_x^{-1})\|_2,
  \eq
  \bq
    \|d_{|z}\si\|= \|d_{|y} (f_\si \circ \psi_x^{-1})\|,
    \eq
    and 
    \bq
     (d_{|z}\si)^*= (d_{|y} (f_\si \circ \psi_x^{-1}))^*,
     \eq
 where $y= \psi_x(z)$. 
 Finally, we denote, for every small enough $r>0$, by $B(x,r)\subset W_x$ the ball centered at $x$ and of
 radius $r$ for the flat metric of $V_x$ pulled back by $\psi_x$, so that 
 \bq \label{ball}
 B(x,r) = \psi_x^{-1}(B(0,r)).
 \eq
\begin{Proposition}\label{prop 3}
Let $X$ be a smooth real projective manifold of dimension $n$, $(L,h_L)$ be
a real holomorphic Hermitian line bundle of positive curvature over $X$ and $(E,h_E)$
be a rank $k$ real holomorphic Hermitian vector bundle, with $1\leq k \leq n$. 
We choose a field of ${h_L}$-trivializations on $\R X$.  Then, for every regular pair $(U,P)$, every  large enough integer  $d$,
 every  $x$ in $\R X$ and every local trivialization of $E$ orthonormal at $x$, 
  there exist $\sigma_\up\in \rhxed $   and an open subset $U_d$ of $ B(x,\frac{R_\up}{\sqrt d})\cap \R X$
such that 
\begin{enumerate}
\item $ \|\si_\up\|_{L^2(h_L)}$ be equivalent to $\frac{\|P\|_{L^2} }{\sqrt{\delta_L}}$ as $d$ grows to infinity,
where $\| P\|_{L^2}$ is defined by (\ref{polyL2}) and $\delta_L = \int_X \omega^n $,  \label{un}
\item $ (U_d, \sigma_\up^{-1} (0) \cap U_d) $ be diffeomorphic  to $(U, P^{-1}(0)\cap U) \subset \R^n$,
\item for every $(\delta,\epsilon)\in \mathcal T_\up$ given by Definition \ref{def 2},  there exists a compact subset $K_d\subset U_d$
such that $$\inf_{U_d\setminus K_d} |\si_\up| >\frac{\delta}{2}\sqrt d^n, $$
while for every $y$ in $U_d$,  
\begin{equation}\label{donaldson}
 |\si_\up (y)|< \frac{\delta}{2} \sqrt d^n \Rightarrow  
\forall w\in \R^k, \ |(d_{|y}\si_\up)^*(w) |\geq \frac{\epsilon}{2}\sqrt d^{n+1} |w|.
\end{equation}
\end{enumerate}
\end{Proposition}
\bpr
We proceed as in the proof of Proposition 3.2 of \cite{GaWe4}. Let $(U,P)$ be a regular pair,
$x\in \R X $ and $d$ large enough. We set $U_d = \psi_x^{-1}(\frac{1}{\sqrt d}U ) \subset B(x,\frac{R_{(U,P)}}{\sqrt d})$
and $K_d = \psi_x^{-1}(\frac{1}{\sqrt d} K ). $ Let $\chi : \C^n \to [0,1]$ be a smooth function with compact
support in $B(0, R_{(U,P)})$, which equals one in a neighbourhood of the origin. Then, let $\si $ be 
the global smooth section of $E\otimes L^d$ defined by $\si_{|X\setminus W_x }= 0$ and
$$ \si_{|W_x} = (\chi \circ \psi_x ) \big(\sum_{j=1}^k P_j(\sqrt d \psi_x) e_j\otimes e^d \big),$$
where $P=(P_1, \cdots, P_k)$ is now considered as a function $\C^n \to \C^k$.
From the $L^2$-estimates of H\"ormander, see \cite{Hormander1} or \cite{MaMa}, 
there exists a global section $\tau$ of $E\otimes L^d$ such that $\dbar \tau = \dbar \si$ and
$\|\tau\|_{L^2(\held)}\leq \| \dbar \si\|_{L^2 (\held)}$ for $d$ large enough. 
This section $\tau$ can be chosen orthogonal to holomorphic sections  and is then unique, in particular real. Moreover, there
exist positive constants $c_1$ and $c_2$, which do not depend on $x$, such that $\| \tau\|_{L^2 (\held)} \leq c_1 e^{-c_2 d}$
and $\sup_{U_d}(|\tau| + \|\tau\|_2)\leq c_2 e^{-c_2d} $, see Lemma 3.3 of \cite{GaWe4}. 
We then set $\si_{(U,P)} = \sqrt d^n (\si - \tau)$. It has the desired properties as can be
checked along the same lines as in the proof of Proposition 3.2 of \cite{GaWe4} 
and thanks to Lemma \ref{lemme 5}.\epr
\begin{Lemma}\label{lemme 5}
Let $U$ be an open subset of $\R^n $, $1\leq  k \leq n$, $f: U\to \R^k $  be a  function of class $C^1$  and $(\delta, \epsilon)\in (\R^*_+)^2$
be such that
\begin{enumerate}
\item there exists a compact subset $K$ of $ U$ such that $\inf_{U\setminus K} |f| >\delta$,
\item for every $ y$ in $U$, $|f(y)|< \delta \Rightarrow 
\forall w\in \R^k, \ |(d_{|y}f)^*(w) |\geq \epsilon  |w|. $
\end{enumerate}
Then, for every function $g : U\to \R^k $ of class $C^1$ such that $\sup_U|g| <\delta$ and $\sup_U \|dg\| <\epsilon$,
zero is a regular value of $f+g$ and $(f+g)^{-1}(0)$ is compact and isotopic to $f^{-1}(0)$ in $U$.
\end{Lemma}
\bpr
The proof is analogous to the one of Lemma 3.4 of \cite{GaWe4}, since $ \| (dg)^* \| = \| dg\|.$
\epr

The following Lemma \ref{lemma tian} establishes the existence of peak sections for
higher rank vector bundles. 
\begin{Lemma}[compare Lemma 1.2 of \cite{Tian}]\label{lemma tian}
Let $X$ be a smooth real projective manifold of dimension $n$, $(L,h_L)$ be
a real holomorphic Hermitian line bundle of positive curvature over $X$ and $(E,h_E)$
be a rank $k$ real holomorphic Hermitian vector bundle, with $1\leq k \leq n$.
Let $x\in \R X$, $(p_1, \cdots, p_n)\in \Nn^n$, $i\in\{1, \cdots, k\}$
and $p' >p_1+\cdots + p_n$. There exists $d_0\in \Nn$ independent of $x$ such that for every $ d>d_0$,
there exists $\si\in \rhxed $ with the property that $\|\si\|_{L^2(h_L)}= 1$ and
 if $ (y_1, \cdots, y_n)$ are  local real holomorphic coordinates in the neighbourhood 
of $x$ and $(e_1, \cdots e_k)$ is a local real holomorphic trivialization of $E$
orthonormal at $x$,  we can assume that in a neighbourhood of $x$,
\begin{equation}\label{sigma}
 \si(y_1, \cdots, y_n) = \lambda y_1^{p_1}\cdots y_n^{p_n} e_i\otimes e^d (1+ O(d^{-2p'} )) + O(\lambda |y|^{2p'}), 
 \end{equation}
where 
$ \lambda^{-2}= \int_{B(x,\frac{\log d }{\sqrt d})}|y_1^{p_1}\cdots y_n^{p_n}|^2
\|e^d\|^2_{h_L^d} dV_{h_L}$, 
with $dV_{h_L}= \omega^n/{\int_X \omega^n}$
 and where $e$ is a local trivialization of $L$ 
whose potential $ -\log h_L (e,e)$ reaches a local minimum at $x$ with
Hessian $\pi \omega(.,i.)$.
\end{Lemma} 
\bpr 
The proof goes along the same lines as the one of Lemma 1.2 of \cite{Tian}. 
Let  $\eta$ be a cut-off function on $\R$ with $\eta = 1$ in a neighbourhood of $0$,  
and $$\psi = (n+2p') \eta \big(\frac {d\|z\|^2}{\log^2 d }\big) \log \big(\frac {d\|z\|^2}{\log^2 d}\big) $$
 in the coordinates $z$ on $X$. 
Then, $i\partial \dbar  \psi$ is bounded from below by $-C\omega$, where $C$ is
 some uniform constant independent of $d$ and $x$. 
Let $s\in C^{\infty}(X, E\otimes L^d)$ be the real section defined by 
$$s = \eta \big(\frac {d\|z\|^2}{\log^2 d }\big) y_1^{p_1}\cdots y_n^{p_n} e_i\otimes e^d. $$ 
Then,  from Theorem 5.1 of \cite{Demailly82}, 
for $d$ large enough not depending on $x$,
there exists a real section $u\in C^\infty (X, E\otimes L^{d})$ 
such that $\dbar u = \dbar s$ and satisfying 
the H\"ormander  $L^2$-estimates
$$ \int_X \|u\|_{\held}^2 e^{-\psi} dV_{h_L} \leq \int_X \|\dbar s \|_{\held}^2 e^{-\psi} dV_{h_L}.$$
The presence of the singular weight $e^{-\psi}$ forces the jets of $u$ to vanish up to order $2p'$ at $x$. 
As in Lemma 1.2 of \cite{Tian}, we conclude that the real holomorphic section $\si = {(s-u)}/{\|s-u\|_{L^2(\held)}}$
satisfies the required properties.
\epr 
In this first section we will only need peak sections given by Lemma \ref{lemma tian} with $\sum_{i=1}^n p_i=0$, whereas
in the second one we will need those given with $\sum_{i=1}^n p_i\leq 2$. 
\begin{Definition}\label{sections pics}
For $i\in \{1, \cdots, k\}, $ let  $\si^i_0$ be the section given by  Lemma \ref{lemma tian} 
with $p'=3$ and $p_1 = \cdots = p_n = 0$. Likewise, for every $j\in \{1, \cdots, n\}$, let $\si^i_j$ be a section
given by (\ref{sigma})  with $p'=3$, $p_j= 1$ and $p_l= 0$ for $l\in \{1, \cdots, n\}\setminus \{j\}$.
Finally, for every $1\leq l\leq m\leq n$, let $\si^i_{lm}$ be a section given by (\ref{sigma})
with $p'=3$, $p_j= 0$ for every $j\in \{1, \cdots , n\} \setminus \{l,m\}$ and $p_l = p_m = 1$ 
if $l\neq m$, while $p_l=2$ otherwise. 
\end{Definition}
The asymptotic values of the constants $\lambda$ in (\ref{sigma}) are
given by Lemma \ref{LemmaTian2} (compare Lemma 2.1 of  \cite{Tian}).
\begin{Lemma}\label{LemmaTian2}
For every $i\in \{1, \cdots, k\}$, the sections given by Definition \ref{sections pics} satisfy
\begin{eqnarray}
\si^i_0/\sqrt{\delta_L d^n} & \equid & e_i\otimes e^d  + O(\|y\|^6), \label{sizero}\\
\forall j\in \{1, \cdots, n\}, \si^i_j/\sqrt{\pi \delta_L d^{n+1} } & \equid  &  y_j e_i\otimes  e^d  + O(\|y\|^6),  \label{sigmaj}\\
\forall l,m\in \{1, \cdots, n\}, l\neq m, \si^i_{lm}/\big(\pi \sqrt{\delta_L d^{n+2}}\big)&\equid  & y_ly_m e_i\otimes e^d  + O( \|y\|^6 ), \label{sigmakl} \\
\text{and }\forall l\in \{1, \cdots, n\}, \si^i_{ll}/\big(\pi \sqrt{\delta_L d^{n+2}}\big) &\equid & \frac{1}{\sqrt 2} y_l^2e_i\otimes e^d  + O(\|y\|^6).\label{sigmakk}
\end{eqnarray}
\end{Lemma}
Moreover, these sections
are asymptotically orthonormal as $d$ grows to infinity, as follows from Lemma \ref{LemmaTian3}.
\begin{Lemma}[compare Lemma 3.1 of \cite{Tian}]\label{LemmaTian3}
For every $x\in \R X$, 
the sections $(\si^i_j)_{\underset{ \ 0\leq j \leq n}{1\leq i \leq k}}$ and $(\si^i_{lm})_{\underset{1\leq l\leq m \leq n}{1\leq i\leq k}}$ 
given by Definition \ref{sections pics}
have $L^2$-norm equal to one
and their pairwise scalar product are dominated by a $O(d^{-1})$ which does not depend on $x$.
Likewise, their scalar products with every section of $\rhxed$ of $L^2$-norm equal to one and whose 2-jet at $x$ vanishes is 
dominated by a 
 $O(d^{-3/2})$ which does not depend on $x$. 
\end{Lemma}
\bpr The proof goes along the same lines as the one of Lemma 3.1 of \cite{Tian}.
\epr
\begin{Lemma}\label{other}
Denote by  $v$ the density of $dV_{h_L}= \omega^n /\int_X \omega^n $ with respect to the 
 volume form $dx$ chosen in  (\ref{l2product}), so that 
$dV_{h_L} = v(x)dx$. Then the sections given by Defintion \ref{sections pics} times $\sqrt {v(x)}$
are still asymptotically orthonormal for (\ref{l2product}).
\end{Lemma}
\bpr This is a direct consequence of Lemmas \ref{lemma tian} and \ref{LemmaTian3} and the 
asymptotic concentration of the support of the peak sections near $x$. \epr
\begin{Remark}\label{rkcx}
The complex analogues of Lemmas \ref{lemma tian}, \ref{LemmaTian2} and \ref{LemmaTian3} 
hold, compare \cite{Tian}.
\end{Remark}

\subsection{Proof of Theorem \ref{th2}}
 We first compute the expected local $C^1$-norm of sections.
\begin{Proposition}\label{prop 4}
Let $X$ be a smooth real projective manifold of dimension $n$, $(L,h_L)$ be
a real holomorphic Hermitian line bundle of positive curvature over $X$ and $(E,h_E)$
be a rank $k$ real holomorphic Hermitian vector bundle, with $1\leq k \leq n$. 
We equip 
$\R X$ with a field of ${h_L}$-trivializations, see \S \ref {para II.1}. Then, for every positive  $R$,
\beq
\limsup_{d\to \infty} \sup_{x\in \R X}\frac{1}{d^n} E\big(\sup_{B(x,\frac{R}{\sqrt d})} \frac{|\si|^2}{v(x)}\big) 
& \leq & 6 k  \delta_L  \rho_R  \text{ and}\\
\limsup_{d\to \infty} \sup_{x\in \R X} \frac{1}{d^{n+1}} E\big(\sup_{B(x,\frac{R}{\sqrt d})} \frac{\|d\si\|_2^2}{v(x)}\big) 
& \leq & 6 \pi n k  \delta_L  \rho_R,
\eeq
where $v$ is given by Lemma \ref{other} and $\rho_R$ is given by (\ref{rhoR}), see (\ref{bal}) and (\ref{bal2})
for the definitions of $|\si|$ and $\|d\si\|_2$. 
\end{Proposition}
Note that a global estimate  on the sup norm of $L^2$ random holomorphic sections is given by Theorem $1.1$ of \cite{SZ4}.

\bpr
The proof goes along the same lines as the proof of Proposition 3.5 of \cite{GaWe4}. 
We first establish from the mean value inequality that 
for every $x\in \R X$,  $R>0$ and $s>0$,
$$E \big(\sup_{B(x,\frac{R}{\sqrt d})}|\si|^2 \big)\leq \frac{1}{Vol(B(\frac{s}{\sqrt d}))}\int_{B(x,\frac{R+s}{\sqrt d})} \E( |\si|^2 ) \psi_x^*dy$$
for $d$ large enough not depending on $x$. 
Then, for every $z \in B(x, \frac{R+s}{\sqrt d})\cap \R  X$, we write $\si = \sum_{i=1}^k a_i \si_0^i + \tau$,
where $\tau \in \rhxed$ vanishes at $z$ and $(\si_0^i)_{i= 1, \cdots k}$ are the peak  sections 
at $z$  given by Definition \ref{sections pics}. 
In particular, by Lemma \ref{LemmaTian2}, at the point $z$, for every $i=1, \cdots, k$, 
$\|\si_0^i\|_{h_{E,d}} \equid \sqrt{\delta_L d^n}$. 
Moreover, since $(e_1, \cdots, e_n)$ is orthonormal at $x$,
\beq | \si_0^i (z)|^2 &=& \|\si_0^i (z)\|^2_{h_{E,d}} (1+ O(|z-x|) e^{d\phi(z)} \\
& \leq & \delta_L d^n e^{\pi (R+s)^2} (1+o(1))
\eeq
from the inequalities (\ref{phi}), 
where the $o(d^n)$ can be chosen not to depend on $x\in \R X$. 
Suppose that $dy = dV_{h_L}$. Then, 
by Lemma \ref{LemmaTian3}, the peak sections are
asymptotically orthogonal to each other for the scalar product defined by (\ref{l2product}), and asymptotically orthogonal
to the space of sections $\tau $ vanishing at $x$. 
We deduce that
\beq \E(|\si(z)|^2)
& = & \E(\big | \sum_{i=1}^k  a_i \si_0^i\big|^2) \, (1+ o(1))\\
& = & (\sum_{i=1}^k |\si_0^i(z)|^2 )\frac{1}{\sqrt \pi} \int_\R a^2 e^{-a^2}da \, (1+ o(1))\\
& \leq & \frac{1}{2} k \delta_L d^n e^{\pi (R+s)^2 }(1+ o(1)).
\eeq
When $z\notin B(x,\frac{R+s}{\sqrt d})\cap \R X$, the space of real sections
vanishing at $z$ gets of real codimension $2k$ in $\rhxed$. 
Let $\langle \theta^i_1, \theta^i_2, \ i\in \{1, \cdots, k\} \rangle $ be
an orthonormal basis of its orthogonal complement. From Remark \ref{rkcx},
for every $i\in \{1, \cdots, k\}$,  $j\in \{1,2\}$,
$$\limsup_{d\to \infty} \frac{1}{d^n } |\theta_j^i(z)|^2 \leq 2 \delta_L e^{\pi(R+s)^2},$$
an upper bound which does not depend on $z$.
We deduce that 
\beq \E (|\si(z)|^2) &=& \int_{\R^{2k} }|\sum_{i=1}^k (a^i_{01} \theta_1^i (z) + a^i_{02} \theta^i_2(z))|^2 e^{-\sum_{i=1}^k (a_{01}^i)^2+ 
(a_{02}^i)^2} \frac{1}{\pi^k} \Pi_{i=1}^k da^i_{01} da^i_{02} \\
& \leq & 2\delta_L d^n e^{\pi (R+s)^2}(1+o(1))\sum_{i=1}^k \int_{\R^2} \big((a^i_{01})^2+ (a^i_{02})^2 + 2|a^i_{01}| |a^i_{02}| \big) \dots\\
&& \dots \frac{1}{\pi}e^{-(a_{01}^i)^2 -(a_{02}^i)^2} da^i_{01}   da^i_{02}\\
& \leq & 6 \delta_L d^n e^{\pi (R+s)^2}(1+o(1)).
\eeq

We deduce the first part of Proposition \ref{prop 4} by taking the supremum over $\R  X$, choosing $s$ 
which minimize $g_{R_{(U,P)}}$ and taking the $\limsup$ as $d$ grows to infinity.

 In general, the Bergman section  at $x$ for the $L^2$-product (\ref{l2product}) 
associated to the volume form $dx$
is equivalent to the Bergman section $\si_0$ at $x$ for $dV_h$ times $\sqrt {v(x)}$,
see Lemma \ref{other}.
The same  holds true for 
the $\sigma_j$'s, and  the result follows by replacing $\delta_L$ with $v(x)\delta_L $.

The proof of the second assertion goes along the same lines, see the proof of Proposition 3.5 of \cite{GaWe4} (and \cite{SZ5} for similar
results).
\epr
As in \cite{GaWe4}, we then compute the probability of presence of closed affine real algebraic submanifolds,
inspired by an approach of Nazarov and Sodin \cite{Nazarov-Sodin}, see also \cite{Lerario-Lundberg}.
Let $(U,P)$ be a regular pair given by Definition \ref{def 1} and $\Sigma = P^{-1}(0)\subset U$.
Then, for every $x\in \R X $, we set $B_d = B(x,\frac{R_{(U,P)}}{\sqrt d })\cap \R X$, see (\ref{ball}), and denote by 
$Prob_{x,\Si}(E\otimes L^d)$ the probability that $\si \in \rhxed $ has the property that 
$\si^{-1}(0)\cap B_d$ contains a closed submanifold $\Sigma'$ 
such that the pair $(B_d, \Sigma')$ be diffeomorphic to $(\R^n, \Si)$. That is,
\beq
  Prob_{x,\Sigma}(E\otimes L^d) = \mu_\R \left\{\si \in \rhxed \, | \,
( \si^{-1}(0)\cap B_d) \supset \Sigma', \ (B_d,\Sigma') \sim (\R^n,\Sigma)  \right\}.
 \eeq
We then set $Prob_{\Sigma}(E\otimes L^d) = \inf_{x\in \R X} Prob_{x,\Sigma}(E\otimes L^d)$.

\begin{Proposition}\label{prop 5}
Let $X$ be a smooth real projective manifold of dimension $n$, $(L,h_L)$ be
a real holomorphic Hermitian line bundle of positive curvature over $X$ and $(E,h_E)$
be a rank $k$ real holomorphic Hermitian vector bundle, with $1\leq k \leq n$. Let $(U,P)$ be
a regular pair given by Definition \ref{def 1} and $\Sigma= P^{-1}(0)\subset U.$ 
Then, 
$$\liminf_{d\to \infty} Prob_{\Sigma}(E\otimes L^d) \geq m_{\tau_\up}, $$  see (\ref{mtau}).
\end{Proposition}
\bpr
The proof is the same as the one of Proposition 3.6 of \cite{GaWe4} and is not reproduced here.
\epr
The proof of Theorem \ref{th2} (resp. Corollary \ref{coro 1}) then just goes along
the same lines as the one of Theorem 1.2 (resp. Corollary 1.3) of \cite{GaWe4}.

\subsection{Proof of Theorem \ref{theo 0}}
 Let $\up$ be a regular pair given by Definition \ref{def 1}. For every $d>0$,
let $\Lambda_d$ be a  maximal  subset of $\R X$ with the property that two distinct points of $\Lambda_d$
are at distance greater than $\frac{2R_\up}{\sqrt d}$. The balls centered at points of $\Lambda_d$
and of radius $\frac{R_\up}{\sqrt d}$ are disjoints, whereas the ones of radius
$\frac{2R_\up}{\sqrt d}$ cover $\R X$. 
Note that if we use the local flat metric given by a trivial $h_L$-trivialization, 
then the associated lattice has  asymptotically the same number of balls
than $\Lambda_d$ as $d$ grows to infinity, so we can suppose from now on that the balls 
are defined for this local metric. 
For every $\si \in \rhxed$, denote by $N_\Sigma(\Lambda_d, \sigma)$ 
the number of $x\in \Lambda_d$ such that 
 the ball $B_d = B(x,\frac{R_\up}{\sqrt d})\cap \R  X$ contains a codimension $k$ submanifold
$\Sigma'$ with $ \Sigma' \subset \si^{-1}(0)$ and $(B_d,  \Sigma')$
 diffeomorphic to $(\R^n, \Sigma)$.
By definition of 
$N_\Sigma(\si) $, 
$N_\Sigma (\Lambda_d, \sigma) \leq N_\Sigma(\sigma),$
see \S \ref{intro lower},
while from Proposition \ref{prop 5}, for every $0<\epsilon <1$,
\beq
 |\Lambda_d|   m_{\tau_\up}  &\leq &\sum_{x\in \Lambda_d} Prob_{x,\Sigma}(E\otimes L^d)\\
&\leq &  \sum_{j=1}^{ |\Lambda_d|}  j
 \mu_\R \{\si |  N_{\Sigma} (\Lambda_d, \si)= j\}  \\
 & \leq & \epsilon m_{\tau_\up}  |\Lambda_d| \mu_\R 
 \left\{\si | N_{\Sigma} (\Lambda_d, \si) \leq \epsilon m_{\tau_\up}  |\Lambda_d|\right\} \\
& &  + |\Lambda_d| \mu_\R 
  \left\{\si | N_{\Sigma} (\Lambda_d, \si) \geq \epsilon m_{\tau_\up}  |\Lambda_d|\right\}.
 \eeq
 We deduce that 
\bq\label{minoration} 
(1-\epsilon )m_{\tau_\up} \leq \mu_\R \left\{ \si|\, N_\Sigma(\si )\geq \epsilon m_{\tau_\up} |\Lambda_d|\right\}
\eq
 and the result follows by choosing a sequence $(U_p,P_p)_p \in \mathcal I_\Sigma$
  such that $$\lim_{p\to \infty} m_{\tau_{(U_p,P_p)}}|\Lambda_d| = c_\Sigma Vol_{h_L} (\R X) \sqrt d^n ,$$
 see (\ref{csigma}).
$\Box$
\subsection{Proof of Corollary \ref{coro i}}\label{para exemples}
In this paragraph, for every positive integer $p$, 
$S^p$ denotes the unit sphere in $\R^{p+1}.$ 
Corollary \ref{coro i} is a consequence of Theorem \ref{th2} and
the following Propositions \ref{prop 1} and \ref{prop i}.
\begin{Proposition}\label{prop 1}
For every $1\leq k\leq n$, $c_{S^{n-k}}\geq \exp (-e^{54+5n}  ).$
\end{Proposition}
Recall the following.
\begin{Lemma}[Lemma 2.2 of \cite{GaWe4}]\label{lemma poly} 
If $P = \sum_{(i_1, \cdots, i_n)\in \Nn^n} a_{i_1, \cdots, i_n} z^{i_1}_1\cdots z_n^{i_n} \in \R [z_1, \cdots, z_n]$,
then
$$ \| P\|^2_{L^2} = \int_{\C^n} |P(z)|^2 e^{-\pi |z|^2} dz = 
\sum_{(i_1, \cdots, i_n)\in \Nn^n} |a_{i_1, \cdots, i_n} |^2 \frac{i_1!\cdots i_n !}{\pi^{i_1+\cdots +i_n}}.$$
\end{Lemma}

\noindent
\bpr[ of Proposition \ref{prop 1}]
For every $n>0$, we set $P_k(x_1, \cdots, x_n) = \sum_{j=k}^n x_j^2 - 1$. 
 For every $x\in \R^n$ and $\delta>0$,
$$
|P_k(x)| <\delta \Leftrightarrow 1 -\delta < \sum_{i=k}^n x_i^2<  1+\delta 
 \Rightarrow  \|d_{|x} P_k\|_2^2 = 4 \sum_{i=k}^n x_i^2 > 4(1-\delta).$$
Moreover from Lemma \ref{lemma poly}, 
$$\|P_k\|^2_{L^2} =1 + \frac{2(n-k+1)}{\pi^2} \leq n-k+2.$$
Now set $P_S = (P_1, \cdots, P_k)$ 
with $P_j(x)= x_j$ for $1\leq j\leq k-1$, so that 
$$\|P_S\|_{L_2}^2 \leq  (k-1)/\pi + (n-k+2) \leq n+1\leq 2n.$$
Since for every $w= (w_1, \cdots , w_k)\in \R^k$ and every $ x\in \R^n $,
$$
 | d_{|x}P_S^*(w)|^2 
 =  \sum_{i=1}^{k-1} w_i^2 + w_k^2 \|d_{|x}P_k\|_2^2,
 $$
we get that $\|d_{|x} P_S^*\|^2\geq \min\big(1, 4(1-\delta)\big)$ if $|P_k(x)|< \delta$.
Choose $$ U_S = \{(x_1, \cdots , x_n)\in \R^n \, | \, \sum_{j=1}^n x_j^2 < 4\}.$$
 Then
if $0<\delta <1$, $$K_\delta = \left\{x\in U_S \ | \ 1 - \delta\leq \sum_{i=k}^n x_i^2\leq 1+\delta \text{ and } \sum_{i=1}^{k-1}x^2_k \leq
1- \frac{1}{2}(1+\delta)^2\right\}$$
is compact in $U_S$ and taking $R^2_\ups = 4$, we see that 
the pair $\ups$ is regular in the sense of Definition \ref{def 1}. The submanifold 
$P_S^{-1}(0) \subset U_S$
is isotopic in $\R^n$ to the unit sphere $S^{n-k}$. 
We deduce  that $(3/4,1)\in \mathcal T_\ups$. 
From (\ref{tauUP}) and (\ref{estimation r})
we deduce
\beq
\tau_\ups  \leq 24 k 4^n e^{16\pi } 2n (2+\pi n)
 \leq  e^{53 + 5n}.
\eeq
The estimate
$c_{S^{n-1}}  \geq \exp (-e^{54+5n}) $ follows then from (\ref{estimation c}).
\epr
\begin{Proposition}\label{prop i}
For every $1\leq k\leq n$ and every $0\leq i\leq n-k$, $c_{S^i \times S^{n-i-k}}\geq \exp (-e^{82+6n}  ).$
\end{Proposition}
\bpr
For every $1\leq k \leq n$ and every $0\leq i\leq n-k$, we set 
 $$Q_k((x_1, \cdots, x_{i+1}),  (y_1, \cdots, y_{n-i-1})) = \big(|x|^2 - 2\big)^2 + \sum_{j=1}^{n-k-i}y_j^2 -1.$$
 For every $(x,y)\in \R^{i+1}\times \R^{n-i-1}$ and $0<\delta<1/2$,
\beq
|Q_k(x,y)| <\delta &\Leftrightarrow& 1 -\delta < (|x|^2-2)^2 + \sum_{j=1}^{n-k-i}y_j^2 <  1+\delta \\
 &\Rightarrow  &\|d_{|(x,y)} Q_k\|_2^2 = 4 \sum_{j=1}^{n-k-i}y_j^2 + 16 |x|^2 (|x|^2 -2)^2,
 \eeq
 with $|x|^2 > 2 - \sqrt {1+\delta} >1/2.$ Thus $\|d_{|(x,y)} Q_k\|_2^2 > 4(1-\delta),$
 compare Lemma 2.6 of \cite{GaWe4}.
Moreover from Lemma \ref{lemma poly}, 
$\|Q_k\|^2_{L^2} \leq 13n^2,$ compare \S 2.3.2 of \cite{GaWe4}.
Now set $Q = (Q_1, \cdots, Q_k)$ 
with $Q_j(x,y)= y_{n-i-j}$ for $1\leq j\leq k-1$, so that 
$$\|Q\|_{L_2}^2 \leq  (k-1)/\pi + 13 n^2\leq 13(n+1)^2.$$
For every $w= (w_1, \cdots , w_k)\in \R^k$ and every $ (x,y)\in \R^{i+1}\times \R^{n-i-1} $,
\beq
 | d_{|(x,y)}Q^*(w)|^2 
 &=&  \sum_{i=1}^{k-1} w_i^2 + w_k^2 \|d_{|(x,y)}Q_k\|_2^2\\
 & > & \min(1,4(1-\delta)) |w|^2 
\eeq
if $|Q_k(x,y)|\leq \delta < 1/2.$
We choose 
$$ U = \{(x,y)\in  \R^{i+1}\times \R^{n-i-1} \, | \, |x|^2 + |y|^2 <6\},$$
$$ K_\delta =\left\{(x,y)\in  U \, | \, 1 - \delta\leq (|x|^2 -2)^2 + \sum_{j=1}^{n-k-i}y_j^2 \leq 1+\delta
\text{ and } \sum_{j=1}^{k-1}y_{n-i-j}^2 \leq 1-\delta \right\},$$
and $R^2_{(U,Q)}= 6$. The pair $(U,Q)$ is regular in the sense of Definition \ref{def 1}
and $Q^{-1}(0)\subset U$ is isotopic in $\R^n$ to the product $S^i \times S^{n-i-k}$
of unit spheres in $\R^{i+1}$ and $\R^{n-i-k+1}$. 
We deduce  that for every positive  $\ep $, $(1/2-\ep,1)\in \mathcal T_{(U,Q)}$ 
and from (\ref{tauUP}) and (\ref{estimation r})
that
\beq
\tau_{(U,Q)}  \leq 24k 4^n e^{24\pi } 13(n+1)^2 (4+\pi n)
 \leq  e^{81 + 6n}.
\eeq
The estimate 
$c_{S^i \times S^{n-i-k}}  \geq \exp (-e^{82+6n}) $ follows then from (\ref{estimation c}).
\epr

\section{Upper estimates for the expected Betti numbers}\label{II}

\subsection{Statement of the results}\label{statement}
For every $1\leq k \leq n$, we denote by $\text{Gr}(k-1,n-1) $ the Grassmann manifold
of $(k-1)$-dimensional linear subspaces of $\R^{n-1}.$ The tangent space of $\grkn$
at every $H\in \grkn$ is canonically isomorphic to the space of linear maps $\lhp$ from $H$ to
its orthogonal $H^\perp$ and we equip it with the norm $$A\in \lhp \mapsto \|A\|_2 = \sqrt {Tr (A^*A)} \in \R^+.$$
The total volume of $\grkn$ for this Riemannian metric is denoted by $Vol(\grkn)$ 
and we set $$\vkn = \frac{1}{\sqrt \pi^{(k-1)(n-k)}}Vol(\grkn) $$ its volume for
the rescaled metric $A\in \lhp \mapsto \frac{1}{\sqrt \pi }\|A\|_2. $ 
Likewise, we equip $M_{k-1}(\R)$ with the Euclidean norm
$A\in \mkr \mapsto \|A\|_2 = \sqrt {Tr (A^*A)}$ and set 
$d\mu(A) = \frac{1}{\pi^{k-1}} e^{-\|A\|_2^2}dA$ 
the associated Gaussian measure on $\mkr$. Then, we set 
$$ E_{k-1}(|\det |^{n-k+2})= \int_{M_{k-1(\R)}} |\det A|^{n-k+2} d\mu(A).$$
\begin{Remark}\label{rk311}
\begin{enumerate}
\item 
The orthogonal group $\text{O}_{n-1}( \R)$ acts transitively on the Grassmannian $\grkn$ with fixators isomorphic to $\text{O}_{k-1}(\R)\times \text{O}_{n-k}(\R).$
We deduce that 
\beq 
Vol\big(\text{Gr}(k-1,n-1)\big) &= &Vol (\text{O}_{n-1}(\R)) / \left(Vol(\text{O}_{k-1}(\R))  \times Vol \big( \text{O}_{n-k} (\R))\right) \\
& = &  {n-1 \choose k-1} \sqrt \pi^{(k-1)(n-k)}\frac{\prod^{k-1}_{j=1}\Gamma(1+ j/2)}{\prod_{j=n-k+1}^{n-1} \Gamma(1+j/2)},
\eeq
where $\Gamma$ denotes the Gamma function of Euler, see for example Lemma 3.4 of \cite{GaWe3}. 
\item From formula (15.4.12) of \cite{Mehta} follows that 
$$E _{k-1}(|\det|^{n-k+2})= \prod_{j=1}^{k-1}\frac{\Gamma(\frac{n-k+2+j}{2})}{\Gamma(\frac{j}{2})},$$
so that $\vkn E_{k-1}(|\det|^{n-k+2})= \frac{(n-1)!}{(n-k)!2^{k-1}}.$
\end{enumerate}
\end{Remark}

We now keep the framework of \S \ref{I.1}. 
 Let us denote, for every $i\in \{0, \cdots, n-k\}$, 
 by $b_i(\rcsi, \R )= \dim H_i(\rcsi, \R)$ the $i$-th Betti number of $\rcsi$ 
 and by $$m_i (\rcsi)= \inf_{f \text{ Morse on } \rcsi}
| \text{Crit}_i (f)| $$ its $i$-th Morse number, where $| \text{Crit}_i (f)| $
denotes the number of critical points of index $i$ of $f$. We then denote by 
$$ \E(b_i) = \int_{\rhxed \setminus \R \Delta_d} b_i(\rcsi, \R) d\mu_\R (\si)$$ 
and $$\E(m_i)= \int_{\rhxed \setminus \R \Delta_d} m_i(\rcsi) d\mu_\R (\si)$$
their expected values. 
The aim of \S \ref{II} is to prove the following Theorem \ref{th1}, see (\ref{ern}) for
the definition of $e_\R (i,n-k-i)$.
\begin{Theorem} \label{th1}
Let $X$ be a smooth real projective manifold of dimension $n$, $(L,h_L)$ be
a real holomorphic Hermitian line bundle of positive curvature over $X$ and $(E,h_E)$
be a rank $k$ real holomorphic Hermitian vector bundle, with $1\leq k \leq n-1$. Then, for every $0\leq i\leq n-k$,
\beq
\limsup_{d\to \infty} \frac{1}{\sqrt d^n } \E(m_i) & \leq &\frac{1}{\Gamma(\frac{k}{2})} \vkn E_{k-1}(|\det|^{n-k+2})
e_\R (i,n-k-i) Vol_{h_L} (\R X).
\eeq 
\end{Theorem}
Note that the case $k=n$ is covered by Theorems \ref{th11} and \ref{theo 3}.
When $k=1$ and $E=\mathcal O_X$, $Vol_{FS} (\R P^k) = \sqrt \pi, $ see Remark 2.14 of \cite{GaWe3}, 
so that Theorem \ref{th1} reduces to Theorem 1.0.1 of \cite{GaWe3} in this case.
The proof of Theorem \ref{th1}
actually goes along the same lines as the one of Theorem 1.1 of \cite{GaWe3}.
The strategy goes as follows. We fix a Morse function $p : \R  X \to \R$. Then, almost
surely on $\si \in \rhxed$, the restriction of $p$ to $\R C_\si$ is itself a Morse function. 
For $i\in \{0, \cdots, n-k\}$, we denote by $\text{Crit}_i(\prcsi)$ the set of critical points of index $i$ of
this restriction and set 
$$\nu_i (\rcsi) = \frac{1}{\sqrt d^n }\sum_{x\in \text{Crit}_i(\prcsi)}\delta_x$$
if $n>k$ and 
$\nu_0(\rcsi)= \frac{1}{\sqrt d^n } \sum_{x\in \rcsi} \delta_x$ if $k=n$. 
We then set 
$$ \E(\nu_i) = \int_{\rhxed} \nu_i(\rcsi) d\mu_\R (\si)$$
and prove the following equidistribution result (compare Theorem 1.2 of \cite{GaWe3}).
\begin{Theorem}\label{theo 3}
Under the hypotheses of Theorem \ref{th1}, let $p : \R X \to \R$ be 
a Morse function. Then, for every $i\in \{0, \cdots, n-k\}$, the measure
$\E(\nu_i)$ weakly converges to $$\frac{1}{\Gamma(\frac{k}{2})} \vkn E_{k-1}(|\det|^{n-k+2})e_\R (i,n-k-i) dvol_{h_L}$$
as $d$ grows to infinity. 
When $k=n$, $\E(\nu_0)$ converges weakly to $\frac{1}{\sqrt \pi}\Gamma(\frac{n+1}{2})
dvol_{h_L}$. 
\end{Theorem}
In Theorem \ref{theo 3}
$dvol_{h_L}$ denotes the Lebesgue measure of $\R  X$ induced by the K\"ahler metric. 
Theorem \ref{th1} is deduced from Theorem \ref{theo 3} by integration of $1$ over $\R X $. 
The next paragraphs are devoted to the proof of Theorem \ref{theo 3}.
\bpr[ of Theorem \ref{th11}] It follows from Theorem \ref{th1},
the Morse inequalities,
 Remark \ref{rk311} and the computation 
 $Vol_{FS} \R P^n = \sqrt \pi/\Gamma(\frac{n+1}{2})$ (see Remark 2.14 of \cite{GaWe3})
when $ k\leq n-1$  and from Theorem \ref{theo 3}
when $k=n$.
\epr
\subsection{Incidence varieties}
Under the hypotheses of Theorem \ref{theo 3}, we set
$$ \rdp = \{\si \in \rhxed | \, \si \in \R \Delta_d \text{ or } \prcsi \text{ is not Morse} \}$$
and 
$$ \mathcal I_i = \{(\si, x)\in (\rhxed \setminus \rdp)\times (\R  X\setminus Crit(p)) \, | \, x\in \text{Crit}_i (\prcsi)\}.$$
We set 
\begin{eqnarray}
\pi_1 : (\si,x)\in \mathcal I_i & \mapsto & \si \in \rhxed \text{ and } \label{pi12}\\
\pi_2 : (\si,x)\in \mathcal I_i & \mapsto & x \in \R X.\label{pi22}
\end{eqnarray}
Then, for every $(\si_0, x_0)\in ((\rhxed \setminus \rdp)\times (\R  X\setminus Crit(p)))$, 
$\pi_1$ is invertible in a neighbourhood $\R U$ of $\si_0$, defining an evaluation map
at the critical point 
$$ ev_{(\si_0, x_0)} : \si \in \R U \mapsto \pi_2 \circ \pi_1^{-1}(\si ) = x\in \text{Crit}_i(\prcsi) \cap \R V,$$
where $\R V$ denotes a neighbourhood of $x_0$ in $\R X $, compare \S 2.4.2 of \cite{GaWe3}. 
We denote by $d_{|\si_0} ev^\perp_{(\si_0,x_0)} $ the restriction of its differential map $d_{|\si_0}ev_{(\si_0,x_0)}$ at $\si_0$
to the orthogonal complement of $\pi_1(\pi_2^{-1}(x_0))$ in $\rhxed$.
\begin{Proposition}\label{prop II 4}
Under the hypotheses of Theorem \ref{theo 3}, $$\E(\nu_i)= \frac{1}{\sqrt d^n } (\pi_2)_*(\pi_1^* d\mu_\R).$$
Moreover, at every point $x\in \R X \setminus Crit(p),$
$$  (\pi_2)_*(\pi_1^* d\mu_\R)_{|x} = \invpin \int_{\pi_1(\pi_2^{-1}(x))} |\det d_{|\si}ev^\perp_{(\si,x)}|^{-1} d\mu_\R (\si) dvol_{h_L}.$$
\end{Proposition}
\bpr
The proof is the same as in the one of Proposition 2.10 of \cite{GaWe3} and is not reproduced here.
\epr
Fix $x\in \R  X\setminus Crit(p)$. Then $\pi_1(\pi^{-1}_2(x))$ is open in a subspace of $\rhxed$. 
Namely,
\begin{eqnarray}\label{24'}
 \pi_1(\pi_2^{-1}(x)) &=& \big\{\si \in \rhxed \setminus \rdp \ | \ \si(x)= 0 \text{ and } \\
 &&\exists \lambda \in \R (E\otimes L^d  )^*_{|x}, \
\lambda \circ \nabla_{|x} \si= d_{|x} p \big\},
\end{eqnarray}
where $\R ((E\otimes L^d  )^*_{|x})$ is the real part of the fibre $ (E\otimes L^d  )^*_{|x}$.
We deduce a well-defined map 
\begin{eqnarray}\label{p}
\rho_x : \pi_1(\pi_2^{-1}(x)) &\to & \text{Gr}(n-k, \ker d_{|x} p)\times (\R (E\otimes L^d  )^*_{|x} \setminus \{0\})\\
\si & \mapsto & (\ker \nabla_{|x} \si, \lambda).
\end{eqnarray}
For every $\si \in \rhxed \setminus \rdp$, the tangent space of $\pi_1 \big(\pi_2^{-1}(x)\big)$
at $\si $ reads 
\beq
T_\si \pix &=& \big\{\sip\in \rhxed \ |\  \sip (x)= 0 \text{ and } \\
&&\exists \lap \in \R(E\otimes L^d)_{|x}^* \ |\  \lap \circ  \nabla_{|x} \si + \lambda \circ \nabla \sip_{| x}= 0\big\}.
\eeq
Likewise, for every $\lambda \in  \R (E\otimes L^d  )^*_{|x}\setminus \{0\}$, the tangent
space of $\rho_x^{-1}(\text{Gr}(n-k, \ker d_{|x} p)\times \{\lambda \})$ at $\si$ reads
$$ T_\si \rho_x^{-1} (\text{Gr}(n-k, \ker d_{|x} p)\times \{\lambda \}) = 
\{\sip \in \rhxed \ |\  \sip(x)= 0 \text{ and } \lambda \circ \nabla_{|x} \sip = 0\}.$$
Finally, for every $K\in \text{Gr}(n-k, \ker d_{|x} p), $ the tangent space of $\rho_x^{-1}(K, \lambda)$ at $\si$
reads 
$$ T_\si \rho_x^{-1} (K,\lambda) = \{\sip \in \rhxed | \sip(x)= 0, \nax \sip_{|K}=0 \text{ and } \lambda \circ \nabla_{|x} \sip = 0\}.$$
Let us choose local real holomorphic coordinates $(x_1, \cdots, x_n)$ of $X$ near $x$ such that $(\partial/\partial x_1, \cdots, \partial/\partial x_n)$
be orthonormal at $x$, with $d_{|x }p $ being colinear to $dx_1$ and such that
$K = \ker \nabla_{|x} \si =  \langle \partial /\partial x_{k+1}, \cdots, \partial/\partial x_n\rangle$. 
Let us choose a local real holomorphic trivialization $(e_1, \cdots, e_k)$ of $E$ near $x$ 
that is orthonormal at $x$ and be such that
$\ker \lambda_{|x} = \langle e_2\otimes e^d , \cdots, e_k\otimes e^d\rangle_{|x}.$
For $d$ large enough,  we define the following subspaces of $\rhxed$:
\begin{eqnarray}
H_x &=& \langle (\si_0^i )_{1\leq i\leq k} , (\si_j^1)_{k+1\leq j\leq n}\rangle   \label{Hx} \\
H_\lambda &=& \langle  (\si_j^1)_{1\leq j\leq k}\rangle \label{Hl} \\ 
H_K &=& \langle (\si_j^i)_{\underset{k+1\leq j\leq n}{2\leq i \leq k}}\rangle, 
 \label{HK}
\end{eqnarray}
where the  sections
$(\si_0^i)_{1\leq i \leq k}$ and $ (\si^i_j)_{\underset{1\leq j\leq n}{1\leq i \leq k}}$ 
of $\rhxed$ are given by Lemma \ref{lemma tian} and Definition \ref{sections pics}. 

$H_K$ is a complement of $T_\si \rho_x^{-1} (K,\lambda)$ in $T_\si \rho_x^{-1}(\text{Gr}(n-k, \ker d_{|x} p)\times \{\lambda\})$,
$H_\lambda$ is a complement of $T_\si \rho_x^{-1}(\text{Gr}(n-k, \ker d_{|x} p)\times \{\lambda\})$ in $T_\si \pix$
and $H_x$ is a complement of $T_\si \pix $ in $\rhxed$. 
 Then, from Lemmas \ref{LemmaTian3} and \ref{other}, 
up to a uniform rescaling by $\sqrt {v(x)}$, these complements are asymptotically orthogonal and their given basis orthonormal. 
Hence, we can assume from now on that $v=1$. 
\begin{Lemma}\label{indices}
Under the hypotheses of Theorem \ref{theo 3}, let $(\si,x)\in \mathcal I_i$ and 
$\lambda\in \R(E\otimes L^d)_{|x}^*\setminus\{0\}$ 
such that $\lambda \circ \nabla_{|x} \si = \dxp$. Then, 
$\lambda \circ \nabla^2\sigma_{|K_x} = \nabla^2(p_{|\R C_\si})_{|x}$, so that the quadratic form 
$\lambda \circ \nabla^2 \sigma_{|K_x}$
is non-degenerated  of index $i$. 
\end{Lemma}
\bpr
The proof is similar to the one of Lemma 2.9 of \cite{GaWe3}.
\epr

\subsection{Computation of the Jacobian determinants}

\subsubsection{Jacobian determinant of $\rho_x$}\label{4.1}
Under the hypotheses of Theorem \ref{theo 3}, let $(\si, x)\in \mathcal I_i$. We set $(K,\lambda)= \rho_x(\si)$
and denote by $d_{|\si} \rho_x^H$ the restriction of $d_{|\si } \rho_x$ to 
$H_K \oplus H_\lambda$. We then denote by $\det (d_{|\si} \rho_x^H)$ the Jacobian determinant of
$d_{|\si} \rho_x^H$ computed in the given basis of $H_\lambda$ and $H_K$, see
(\ref{Hl}), (\ref{HK}) and in orthonormal basis of $T_K \text{Gr}(n-k, \ker d_{|x} p ) \times \R (E\otimes L^d)^*_{|x}.$
By assumption, the operator $\nabla_{|x}\si $ does not depend on the choice of a connection $\nabla$ on 
$E\otimes L^d$ and is onto. We denote by $\nabla_{|x} \si^\perp $ its restriction to the orthogonal
$K^\perp$ of $K = \ker \nabla_{|x} \si $,
\beq
\nabla_{|x} \si^\perp : K^\perp \to \R (E\otimes L^d)_{|x}.
\eeq
 Likewise, for every $(\sip_K, \sip_\lambda)\in H_K \oplus H_\lambda$,
the operators $\nabla_{|x} \sip_K $ and $\nabla_{|x}\sip_\lambda$ do not depend on 
the choice of a connection $\nabla$ on $E\otimes L^d$. 
Finally, we write at a point $y\in \R X$ near $x$
$$ \si (y)= \sum_{i=1}^k \big(a_0^i \si_0^i + \sum _{j=1}^n a^i_j \si_j^i + \sum_{1\leq l\leq m\leq n} a^i_{lm } \si^i_{lm}\big)(y)+ o(|y|^2),$$
where $(a^i_0)$, $(a^i_j)$ and $(a^i_{lm})$
are real numbers and $(\si^i_0)$, $(\si_j^i)$ and $(\si_{lm}^i)$
are given by Definition \ref{sections pics}. 
From Lemma \ref{LemmaTian2} and (\ref{24'}), we deduce
that $a^i_0 = 0 = a^1_j$ for $1\leq i \leq k$ and $k+1\leq j \leq n$,
and that
\begin{equation}\label{a11}
\|\lambda \| \sqrt{\pi \delta_L }\sqrt d^{n+1} |a^1_1| = \| d_{|x}p\| + o(1),
\end{equation}
where the $o(1)$ term is uniformly bounded over $\R X$.
\begin{Lemma}\label{lemma 2}
Under the hypotheses of Theorem \ref{theo 3}, let $(\si,x)\in \mathcal I_i$ and 
$(K,\lambda)= \rho_x (\si). $ Then,
$d_{|x}\rho_x^H $ writes
\beq
H_K \oplus H_\lambda & \to & T_K \text{Gr}(n-k, \kxp )\times \R \elstar_{|x} \\
(\sip_K, \sip_\lambda) & \mapsto & \big(-(\nabla_{|x} \si^\perp)^{-1}_{|\ker \lambda}  \circ \nabla_{|x} \sip_{K|K}   \  , \
-\lambda \circ \nabla_{|x} \sip_\lambda \circ (\nabla_{|x} \si^\perp)^{-1} \big).
\eeq
Moreover, 
$ | \det \ d_{|\si }\rho_x^H|^{-1} = \frac{  |a^1_1|   }{ \|\lambda \|^k }|\det (a^i_j)_{2\leq i,j\leq k}|^{n-k+1} (1+o(1)),$
where the $o(1)$ term is uniformly bounded over $\R X $. 
\end{Lemma}
\bpr 
Let $(\sip_K, \sip_\lambda)\in H_K \oplus H_\lambda $ and $(\si_s)_{s\in ]-\epsilon, \epsilon[}$ be a path 
of $\pix$ such that $\si_0 = \si $ and $\sip_0 = \sip_K + \sip_\lambda$. Then,
for every $s\in ]-\epsilon, \epsilon[$ and every $v_s \in \ker \nabla_{|x}\si_s$, there exists $\lambda_s \in  \R \elstar_{|x}$
such that 
$$ \left\{
\begin{array}{llll}
\nabla_{|x } \si_s (v_s)&=& 0 &\text{ and } \\
\lambda_s \circ \nabla_{|x} \si_s &=& d_{|x} p.&
\end{array}
\right.
$$
By derivation, we deduce
$$ \left\{
\begin{array}{llll}
	\nabla_{|x } \sip_0 (v_0) + \nabla_{|x} \si (\vpo)&=& 0 &\text{ and } \\
\lapo \circ \nabla_{|x} \si + \lambda \circ \nabla_{|x} \sip_0 &=& 0.&
\end{array}
\right.
$$
By setting $\vp$ the orthogonal projection of $\vpo$ onto $K^\perp$, we deduce that 
$$ \left\{
\begin{array}{llll}
\vp &=& -(\nabla_{|x} \si^\perp)^{-1} \circ \nax \sip_K(v_0) & \text{ and }\\
\lapo &=& - \lambda \circ \nax \sip_\lambda \circ (\nax \si^\perp)^{-1}. &
\end{array}
\right.
$$
 The first part of Lemma \ref{lemma 2} 
follows. Now, recall that $d_{|x } p $ is colinear to $ dx_1$, that $K$ is equipped with the orthonormal basis $(\partial/\partial x_{k+1}, \cdots, \partial/\partial x_n)$, 
$K^\perp$ with the orthonormal basis $(\partial/\partial x_{1}, \cdots, \partial/\partial x_{k})$, 
and that
$\ker \lambda_{|x}$ is spanned by  the orthonormal basis $(e_2, \cdots, e_k)_{|x}$. From Lemma \ref{lemma tian}, 
the map 
$$ \sip_K \in H_K \mapsto \nax \sip_{K|K} \in L(K, \ker \lambda)$$
just dilates the norm by the factor 
$\sqrt{\pi \delta_L d^{n+1}}(1+ o(1))$, where the $o(1)$ term is uniformly bounded over $\R X $.  
Now, since the matrix of the restriction of $\nax \si^\perp $ to $K^\perp \cap \ker d_{|x} p$ in the given
basis of $K^\perp\cap \ker d_{|x} p$ and $\ker \lambda$ equals 
$$\sqrt{\pi \delta_L d^{n+1}} (a^i_j)_{2\leq i,j\leq k}
+ o(\sqrt d^{n+1}),$$ 
where the $o(\sqrt d^{n+1})$ term is uniformly bounded over $\R X $. We
deduce that the Jacobian determinant of the map 
$$ M\in L (K, \ker \lambda)\mapsto (\nax \si^\perp_{|\ker \lambda})^{-1}\circ M \in  L(K, K^\perp \cap \ker d_{|x}p))$$
equals 
$$\big((\sqrt {\pi \delta_L d^{n+1}})^{k-1} |\det (a^i_j)_{2\leq i,j\leq k} |(1+o(1))\big)^{k-n}.$$
The Jacobian determinant of the map 
$$ \sip_K \in H_K \mapsto (\nax \si^\perp)^{-1}_{| \ker \lambda} \circ \nax \sip_{K|K} \in T_K \text{Gr}(n-k, \ker d_{|x} p)$$
thus equals  $|\det (a^i_j)_{2\leq i,j\leq k} |^{k-n}+ o(1),$
where the $o(1)$ is uniformly bounded over $\R X$.
Likewise, from Lemma \ref{lemma tian}, the map
$$ \sip_\lambda \in H_\lambda \mapsto \lambda \circ \nax \sip_\lambda \in (K^\perp)^*$$
just dilates the norm by a factor $\sqrt {\pi \delta_L d^{n+1}} \|\lambda\| + o(\sqrt d^{n+1}),$
where the $o(\sqrt d^{n+1})$ is uniformly bounded over $\R X$,
while the Jacobian determinant of the map
$$ M\in (K^\perp)^* \mapsto M\circ (\nax \si^\perp)^{-1}\in \R \elstar_{|x} $$
equals $(\sqrt {\pi \delta_L } \sqrt d^{n+1} )^{-k} |\det (a^i_j)_{1\leq i,j\leq k}|^{-1}(1+o(1))$
so that the Jacobian determinant of the map 
$$\sip_\lambda \in H_\lambda \mapsto \lambda \circ \nax \sip_\lambda \circ (\nax \si^\perp)^{-1}\in \R \elstar_{|x}$$
equals ${\|\lambda \|^k }|\det (a^i_j)_{1\leq i,j\leq k}|^{-1}+ o(1),$
with a $o(1)$ uniformly bounded over $\R X$.
 As a consequence,
$$ | \det  d_{|\si } \rho_x^H|^{-1} = \|\lambda \|^{-k} |\det (a^i_j)_{2\leq i,j\leq k}|^{n-k+1} |a^1_1|(1+o(1)),$$
with a $o(1)$ uniformly bounded over $\R X$,
since the relation $\lambda \circ \nax \si = d_{|x} p$ implies that $a^1_j$ vanishes for $2\leq j\leq n$. 
\epr
\subsubsection{Jacobian determinant  of the evaluation map}
Again, under the hypotheses of Theorem \ref{theo 3} and for $(\si,x)\in \mathcal I_i$, we set 
for every $y$ in a neighbourhood of $x$,
\beqr \label{sigma decomposition} 
\si(y) = \sum_{i=1}^k \big(a_0^i \si_0^i + \sum _{j=1}^n a^i_j \si_j^i + \sum_{1\leq l\leq m\leq n} a^i_{lm } \si^i_{lm}\big)(y)+ o(|y|^2),
\eeqr
where $a^i_0$, $a^i_j$ and $a^i_{lm} $ are real numbers. We then set, for $1\leq l,m\leq n$, 
$\tilde a^1_{ll }= \sqrt 2 a^1_{ll}$,  $\tilde a^1_{lm} = a^1_{lm}$ if $l<m$ and $\tilde a^1_{lm}= a^1_{ml}$
if $l> m$. 
We denote by $d_{|\si}ev^H_{(\si,x)}$ the restriction of $d_{|\si} ev_{(\si,x)} $ to $H_x$, see (\ref{Hx}) and
by $\det  d_{|\si} ev^H_{(\si,x)} $ its Jacobian determinant computed in the given basis of $H_x$ and orthonormal
basis of $T_x \R X$. 
\begin{Lemma}\label{lemma 3}
Under the hypotheses of Theorem \ref{theo 3}, let $(\si,x)\in \mathcal I_i$. Then,
$$ |\det d_{|\si}ev^H_{(\si,x)}|^{-1} = \sqrt {\pi^n d^n  } |a^1_1|  | \det (a^i_j)_{2\leq i,j\leq k}|
| \det (\tilde a^1_{lm})_{k+1\leq l,m\leq n}| (1+o(1)),$$
where the $o(1)$ term has poles of order at most $n-k$ near the critical points of $p$.
\end{Lemma}
\begin{Remark}\label{pole}
In Lemma \ref{lemma 3}, a function $f$ is said to have a pole of order at most $n-k$ near a point $x$ if
$r^{n-k} f $ is bounded near $x$, where $r$ denotes the distance function to $x$. Such a function thus belongs
to $L^1(\R X, dvol_h)$. 
\end{Remark}
\bpr
We choose a torsion free connection $\nabla^{TX}$ (resp. 
a connection $\nabla^{E\otimes L^d}$) on $\R X\setminus Crit(p)$ (resp. on $E\otimes L^d$) such that $\nabla^{TX} dp = 0.$
They induce a connection on $T^*X \otimes \eld $ which makes it possible to differentiate
twice the elements of $\rhxed$. The tangent space of $\mathcal I_i$ then reads
\beqr 
 T_{(\si,x)}\mathcal I_i  &=& \big\{ (\sip, \xp)  \in \rhxed \times T_x \R X \ | \ \sip (x)+ \nabla_{\xp} \si = 0  \text { and } \\
&&\exists \lap \in \R \elstar_{|x}, \ \lap \circ \nax \si + \lambda \circ \nax \sip + \la \circ \nabla^2_{\xp, .} \si = 0\big\}.\label{TIi}
\eeqr
Recall that $T_x \R X $ is the direct sum $K\oplus K^\perp$, where $K = \ker \nax \si$. 
We write $\xp = (\xp_K, \xp_{K^\perp})$ the coordinates of $\xp$ in this decomposition. 
From the first equation we deduce, keeping the notations of \S \ref{4.1}, that
$ \xp_{K^\perp} = -(\nax \si^\perp)^{-1}(\sip (x)).$ From Lemma \ref{lemma tian}, 
the evaluation map at $x$ 
$$\sip \in \langle (\si_0^i)_{1\leq i \leq k} \rangle \mapsto \sip(x) \in E\otimes L^d_{|x}$$
just dilates the norm by a factor 
$\sqrt {\delta_L d^n} (1+o(1))$, , where the $o(1)$ term is uniformly bounded
over $\R X$, 
while $$|\det (\nax \si^\perp)| = (\spd)^k |\det (a^i_j)_{1\leq i,j\leq k}|(1+o(1)).$$
We deduce by composition that the Jacobian 
of the map 
$$\sip \in \langle (\si_0^i)_{1\leq i \leq k} \rangle \mapsto \xp_{K^\perp} = -(\nax \si^\perp)^{-1}(\sip(x)) $$
equals 
$\big(\sqrt{\pi^k d^k }| \det (a^i_j)_{2\leq i,j\leq k}|| a^1_1|\big)^{-1}(1+o(1)),$
where the $o(1)$ term is uniformly bounded over $\R X$. 
Now, equation (\ref{TIi}) restricted to $K$ reads
$$ \la\circ \nabla_{\xp_K, .}^2 \si_{|K} = - \lambda \circ \nax \sip_{|K}.$$
From Lemma \ref{lemma tian}, the map 
$$ \sip \in \langle (\si_j^1)_{k+1\leq j\leq n}\rangle \mapsto - \la \circ \nax \sip_{|K} \in K^*$$
just dilates the norm by a factor $\| \lambda \| \spd (1+o(1))$, with
$o(1)$ term is uniformly bounded over $\R X$.
Likewise, from Lemma \ref{lemma tian}, 
the Jacobian of the map $\lambda \circ \nabla^2 \si_{|K} : K \to K^*$
equals 
\bq \label{etoile}
(\| \la\| \pi \sqrt {\delta_L  d^{n+2}}  )^{n-k}| \det (\tilde a^1_{lm})_{k+1\leq l,m\leq n}| (1+o(1)).
\eq
Here, the $o(1)$ term is no more uniformly bounded over $\R X$ 
though. Indeed, from Lemma \ref{LemmaTian2} and (\ref{sigma decomposition}), 
$$\lambda \circ \nabla^2 \si_{|K}= a^1_1 (\|\lambda\| 
\sqrt{\pi \delta_L d^{n+1}})(\nabla^{TX} dx_1)
+ \sum_{1\leq l\leq m\leq n} \tilde a^1_{lm}  (\|\lambda\| 
\sqrt{\pi \delta_L d^{n+2}})
dx_l \otimes dx_m,$$
since the relation $\lambda \circ \nabla_{|x} \si = \dxp $ imposes
that $a^1_j$ vanishes for $j>1$. Moreover,
since $dp = \sum_{i=1}^n \alpha_i dx_i$, 
with $\alpha_2 (x)= \cdots = \alpha_n (x)= 0$ 
and $|\alpha_1(x)|= \|\dxp\|$, 
we get that 
$$ 0 = \nabla^{TX}(dp)_{|K}= \alpha_1 (\nabla^{TX} dx_1)_{|K} + \sum_{i=1}^n (d\alpha_i \otimes dx_i)_{|K},
$$
so that $\| \nabla^{TX} dx_{1|K}\| = \frac{1}{\|\dxp\|} \| \sum_{i=1}^n d\alpha_i \otimes dx_i\|$
has a pole of order one at $x$. 
In formula (\ref{etoile}), the $o(1)$ term has thus a pole of order at most $n-k$ near the critical points
of $p$. 

We deduce that the Jacobian determinant of the map
$$ \sip \in \langle (\si_j^1 ) _{k+1\leq j\leq n} \rangle \mapsto
 \xp_K = - (\la \circ \nabla^2 \si_{|K })^{-1} \circ (\la \circ \nax \sip_{|K})\in K$$
equals $(\sqrt {\pi ^{n-k} d^{n-k}} |\det (\tilde a^1_{lm})_{k+1\leq l,m\leq n}|)^{-1} (1+o(1))$, up to sign,
where $o(1)$ term has a pole of order at most $n-k$ near the critical points of $p$.
The result follows. 
\epr

\subsection{Proof of Theorem \ref{theo 3}}
\subsubsection{The case $k<n$}
From Proposition \ref{prop 4} we know that 
$$ \E(\nu_i)= \frac{1}{\sqrt {\pi^n d^n}} \big(\int_{\pix} |\det d_{|\si }ev^\perp_{(\si,x)}|^{-1} d\mu_\R (\si) \big ) dvol_{h_L}.$$
From the coarea formula (see \cite{Federer}), we likewise deduce that 
\beq
 \E(\nu_i)&=& \frac{1}{\sqrt{\pi^n d^n }}
 \big( \int_{\text{Gr}(n-k,\ker \dxp) \times \R \elstar_{|x}\setminus\{0\}}
e^{-(a^1_1)^2} 
\frac{ dK\wedge d\lambda }{ \sqrt \pi^{(n-k)(k-1)+k}}...\\
&&... \int_{\rho_x^{-1}(K,\la)} |\det d_{|\si } ev_{(\si,x)}^\perp |^{-1} |\det d_{|\si} \rho_x^\perp|^{-1} d\mu_\R (\si) \big)
dvol_{h_L},
\eeq
since with the notations (\ref{sigma decomposition}),
 $\si \in \rho_x^{-1}(K,\lambda)$ if and only if $\forall i\in \{1, \cdots, k \}$ and $\forall j\in \{k+1, \cdots, n\}$,
 $a^i_0 = 0 = a^i_j$ while $\forall j\geq 2$, $a^1_j=0$ and
 $|a^1_1| = \frac{\| d_{|x}p\| }{\| \la\| \sqrt {\pi \delta_L} \sqrt d^{n+1}}$. 
 From Lemma \ref{LemmaTian3} and the relation (\ref{a11}), we deduce that for every $x\in \R X \setminus Crit(p)$ and every
 $(K,\la)\in \text{Gr}(n-k, \ker \dxp)\times \R \elstar_{|x} \setminus \{0\}$,
 \beq 
 &&\int_{\rho_x^{-1}(K,\la)} |\det d_{|\si } ev_{(\si,x)}^\perp |^{-1} |\det d_{|\si} \rho_x^\perp|^{-1} d\mu_\R (\si)\\
 &\equid & \int_{\rho_x^{-1}(K,\la)} |\det d_{|\si } ev_{(\si,x)}^H |^{-1} |\det d_{|\si} \rho_x^H|^{-1} d\mu_\R (\si).
 \eeq
 Thus, from Lemmas \ref{lemma 2},  \ref{lemma 3} and \ref{indices}, $\E(\nu_i)$ converges to
 \beq 
 && 
 \int_{M_{k-1}(\R)} |\det (a^i_j)_{2\leq i,j\leq k} |^{n-k+2} d\mu(a^i_j)
 \int_{Sym_\R (i,n-k-i)}  |\det (\tilde a^1_{lm})_{k+1\leq l,m\leq n} | d\mu(\tilde a^1_{lm})...\\
 && ... \int_{\text{Gr}(n-k,\ker \dxp) \times \R \elstar_{|x}\setminus\{0\}}
 \frac{ (a_1^1)^2 e^{-(a^1_1)^2} }{\| \lambda \|^k} \frac{ dK\wedge d\lambda }{ \sqrt \pi^{(n-k)(k-1)+k}},
  \eeq
  where the convergence is dominated by a function in $L^1(\R X, dvol_{h_L})$, see Remark \ref{pole}.
We deduce that $\E(\nu_i)$ gets equivalent to 
$$ \frac{\|\dxp\|^2}{\delta_L d^{n+1} \sqrt \pi^{k+2}} \vkn E_{k-1}(|\det |^{n-k+2}) e_\R (i,n-k-i)   \big(\int_{\R \elstar_{|x}\setminus \{0\}}  
\frac{e^{-(a^1_1)^2 } }{\|\lambda\|^{k+2}} d\la   \big) dvol_{h_L}.
$$
Now, 
\beq
\frac{\|\dxp\|^2 }{\pi \delta_L d^{n+1}} \int_{\R \elstar_{|x}\setminus\{0\}}
\frac{e^{-(a^1_1)^2} }{\|\lambda\|^{k+2}} d\la  &= &
\frac{Vol(S^{k-1}) \| \dxp\|^2 }{\pi \delta_L d^{n+1}}\int_0^{+\infty} 
\frac{e^{-(a^1_1)^2} }{\|\lambda\|^{3}} d\|\la\|\\  & = &
Vol(S^{k-1}) \int_0^{+\infty} e^{-r^2} rdr = \frac{1}{2} Vol(S^{k-1}).
\eeq
Since 
$Vol(S^{k-1}) = \frac{2\sqrt \pi^k}{\Gamma (k/2)}$, we finally
deduce that
$\E(\nu_i) $ weakly converges to 
 $$ \frac{1}{\Gamma(k/2)} \vkn E_{k-1}(|\det |^{n-k+2}) e_\R (i,n-k-i) dvol_{h_L},$$
 where the convergence is dominated by a function in $L^1(\R X, dvol_{h_L})$. \ $\Box$

\subsubsection{The case  $k=n$}
When the rank of $E$ equals the dimension of $X$,
the vanishing locus of a generic section $\si$ of $\rhxed$ is 
a finite set of points. 
We set $\nu = \frac{1}{\sqrt d^n }\sum_{x\in \rcsi} \delta_x,$
and define 
the incidence variety as
$$ \mathcal I = \{(\si, x)\in (\rhxed \setminus \R \Delta_d)\times \R  X \, | \, \si(x)= 0\}.$$
The projections $\pi_1$ and $\pi_2$ are defined by  (\ref {pi12}) and (\ref{pi22}).
As before, for every $(\si_0, x_0)\in (\rhxed\setminus \R \Delta_d ) \times \R  X$, 
$\pi_1$ is invertible in a neighbourhood $\R U$ of $\si_0$, defining an evaluation map
at the critical point 
$$ ev_{(\si_0, x_0)} : \si \in \R U \mapsto \pi_2 \circ \pi_1^{-1}(\si ) = x\in \rcsi \cap \R V,$$
where $\R V$ denotes a neighbourhood of $x_0$ in $\R X $, compare \S 2.4.2 of \cite{GaWe3}. 
We denote by $d_{|\si_0} ev^\perp_{(\si_0,x_0)} $ the restriction of its differential map $d_{|\si_0}ev_{(\si_0,x_0)}$ at $\si_0$
to the orthogonal complement of $\pi_1(\pi_2^{-1}(x_0))$ in $\rhxed$.
Then,
from Proposition \ref{prop II 4},
 $$\E(\nu)=\frac{1}{\sqrt d^n }   (\pi_2)_*(\pi_1^* d\mu_\R)_{|x} = \frac{1}{\sqrt {\pi d}^n} \int_{\pi_1(\pi_2^{-1}(x))} |\det d_{|\si}ev^\perp_{(\si,x)}|^{-1} d\mu_\R (\si) dvol_{h_L}.$$
The space $H_x = \langle (\si_0^i)_{1\leq i \leq k}\rangle$ is a complement to $T_\si \pi_1(\pi_2^{-1}(x))$
in $\rhxed$ and in the decomposition (\ref{sigma decomposition}),
$a_0^i=0$ for every $i=1, \cdots, k$. 
The tangent space of $\mathcal I$ at $(\si,x)$
reads
$$ T_{(\si,x)}\mathcal I= \{(\sip, \xp)\in \rhxed\times T_x \R X \, |  \, \sip (x)+ \nabla_{|x}\si (\xp)= 0  \}.$$
As in the proof of Lemma \ref{lemma 3}, we deduce that the  Jacobian determinant  of the map
$$ \sip \in H_x \mapsto \xp=- (\nabla_{|x} \si^\perp )^{-1}(\sip (x)) \in T_x \R X $$
equals $ \sqrt{\pi^n d^n }|\det (a^i_j)_{1\leq i,j\leq n}|(1+o(1))$, so that 
$$|\det d_{|\si} ev^H_{(\si,x)} |^{-1} = \sqrt{\pi d^n } |\det (a^i_j)_{1\leq i,j\leq n }|
(1+o(1)),$$ where the $o(1)$ 
term is uniformly bounded over $\R X$. 
From lemma \ref{LemmaTian3} we deduce that $\E(\nu)$ 
gets equivalent to 
$$ \big(\int_{M_n(\R)} |\det (a^i_j)_{1\leq i,j\leq n }| d\mu(a^i_j)\big)dvol_{h_L}
= E_n(|\det|) dvol_{h_L}.$$
Formula (15.4.12) of \cite{Mehta}, see Remark \ref{rk311}, now gives 
$$E_n(|\det|)=\frac{\Gamma(\frac{n+1}{2}) }{ \Gamma(1/2)} = 
 \frac{1}{Vol_{FS}(\R P^n) },
$$
see Remark 2.14 of \cite{GaWe3},
hence the result. $\Box$

\subsection{Equidistribution of critical points in the complex case}
Let $X$ be a smooth complex projective manifold of dimension $n$, $(L,h_L)$ be
a  holomorphic Hermitian line bundle of positive curvature $\omega$ over $X$ and $(E,h_E)$
be a rank $k$ holomorphic Hermitian vector bundle, with $1\leq k \leq n$. 
For every $d>0$, we denote by $L^d$ the $d$th tensor power of $L$ and by $h^d$ 
the induced Hermitian metric on $L^d$. We denote by $H^0(X,L^d)$ its complex
vector space of global holomorphic sections and by $N_d$ the dimension of $H^0(X,L^d)$. 
We denote then by $\langle  .,. \rangle$  the $L^2$-Hermitian product on this vector space, 
defined by the relation 
\begin{equation}\label{produit}
\forall \si, \tau \in H^0(X,L^d), \langle \si, \tau \rangle = \int_X h^d(\si, \tau) dx.
\end{equation}
The associated Gaussian measure is denoted by $\mu_\C$. It is defined, for every open subset
$U$ of $H^0(X,L^d)$, by 
\begin{equation}\label{mesureC}
 \mu_\C (U) = \frac{1}{\pi^{N_d}} \int_U e^{-\|\si\|^2} d\si,
 \end{equation}
where $d\si $ denotes the Lebesgue measure of $H^0(X,L^d)$. 
For every $d>0$, we denote by $\Delta^d$ the discriminant hypersurface
of $\hxed$, that is the set of sections $\si \in \hxed$ which 
do not vanish transversally. 
For every $\si \in \hxed\setminus \{0\}$, we denote by $C_\si$
the vanishing locus of $\si $ in $X$. For every $\si \in \hxed \setminus \Delta^d$,
$C_\si$ is then a smooth codimension $k$ complex submanifold of $X$. We equip $X$  with a Lefschetz pencil $p:X \dashrightarrow \cpun$.
We then denote, for every $d>0$, by $\Delta^d_p$ the set of sections $\si \in \hxed$) such that  $\si \in \Delta^d$,
or  $\csi$ intersects the critical locus of $p$, or the restriction of $p$ to
 $C_\si$ is not a Lefschetz pencil. 
For $d$ large enough, this extended discriminant locus is of measure 0 
for the measure $\mu_\C$.

For every $\si \in \hxed\setminus \Delta_p^d$, we denote by $Crit(p_{|C_\si}) $ 
the set of critical points of the restriction of $p$ to $C_\si$ and set, for $1\leq k\leq n-1$,
\bq\label{moyenne nu}
\nu(C_\si) = \frac{1}{ d^n} \sum_{x\in Crit(p_{|\csi} )} \delta_x,
\eq
where $\delta_x$ denotes the Dirac measure of $X$ at the point $x$. When $k=n$, $\nu (C_\si)=  \frac{1}{ d^n} \sum_{x\in \csi} \delta_x$.

\begin{Theorem}\label{lefschetz}
Let $X$ be a smooth complex projective manifold of dimension $n$, $(L,h_L)$ be
a  holomorphic Hermitian line bundle of positive curvature $\omega$ over $X$ and $(E,h_E)$
be a rank $k$ holomorphic Hermitian vector bundle, with $1\leq k \leq n$. 
Let $p : X \dashrightarrow \cpun$ be a Lefschetz pencil. Then, the measure $\E(\nu)$ 
defined by (\ref{moyenne nu}) weakly  converges 
to ${n-1 \choose k-1} \omega^{n}$ as $d$ grows to infinity.
\end{Theorem}
When $k=1$,  Theorem \ref{lefschetz} reduces to
Theorem 3 of \cite{GaWe2}, see 
also Theorem 1.3 of \cite{GaWe3}.
\bpr
The proof goes along the same lines as the one of Theorem \ref{th1},
so we only give a sketch of it. 
Firstly, the analogue of Proposition \ref{prop II 4} 
provides
$$\E(\nu)= \frac{1}{ d^n } (\pi_2)_*(\pi_1^* d\mu_\C),$$
and  at every point $x\in X \setminus (\text{Crit}(p)\cup \text{Base}(p)),$
where $\text{Base}(p)$ denotes the base locus
of $p$,
$$  (\pi_2)_*(\pi_1^* d\mu_\C)_{|x} = \frac{1}{\pi ^n} \int_{\pi_1(\pi_2^{-1}(x))} |\det d_{|\si}ev^\perp_{(\si,x)}|^{-2} d\mu_\R (\si) 
\frac{\omega^n}{n!},$$
see Proposition 2.10 of \cite{GaWe3}.
Choosing complex coefficients in decomposition (\ref{sigma decomposition}),
Lemmas \ref{lemma 2} and \ref{lemma 3} 
remain valid in the complex setting, see Remark \ref{rkcx}.
We deduce that 
\beq
 \E(\nu)&=& \frac{1}{{\pi^n d^n}} \big(\int_{\pix} |\det d_{|\si }ev^\perp_{(\si,x)}|^{-2} d\mu_\C (\si) \big ) \frac{\omega^n}{n!} \\
&\equid& \frac{1}{{\pi^n d^n }}
 \big( \int_{\text{Gr}_\C(n-k,\ker \dxp) \times \elstar_{|x}\setminus\{0\}}
e^{-|a^1_1|^2} 
\frac{ dK\wedge d\lambda }{ \pi^{(n-k)(k-1)+k}}...\\
&&... \int_{\rho_x^{-1}(K,\la)} |\det d_{|\si } ev_{(\si,x)}^\perp |^{-2} |\det d_{|\si} \rho_x^\perp|^{-2} d\mu_\C (\si) \big)
\frac{\omega^n}{n!},
\eeq
with $
|a^1_1|$ given by (\ref{a11}), see Lemma \ref{LemmaTian3}
 as before.
 Here, $\text{Gr}_\C(n-k,\ker \dxp)$ denotes the  Grassmann manifold of $n-k$-dimensional
 complex linear subspaces of $\ker \dxp$. 
 From the complex versions of Lemma \ref{LemmaTian2} and \ref{LemmaTian3}, see Remark \ref{rkcx}
 and the relation (\ref{a11}), we deduce that for every $x\in X \setminus (\text{Crit}(p)\cup \text{Base}(p))$ and every
 $(K,\la)\in \text{Gr}(n-k, \ker \dxp)\times \elstar_{|x} \setminus \{0\}$,
 \beq 
 &&\int_{\rho_x^{-1}(K,\la)} |\det d_{|\si } ev_{(\si,x)}^\perp |^{-2} |\det d_{|\si} \rho_x^\perp|^{-2} d\mu_\C (\si)\\
 &\equid& \frac{|a_1^1|^4 {\pi^n d^n }}{\| \lambda \|^{2k}}
 \int_{M_{k-1}(\C)} |\det (a^i_j)_{2\leq i,j\leq k} |^{2(n-k+2)} d\mu(a^i_j)... \\
&&... \int_{Sym_\C (n-k)}  |\det (\tilde a^1_{lm})_{k+1\leq l,m\leq n} |^2 d\mu(\tilde a^1_{lm}).
  \eeq
We deduce that $\E(\nu)$ is equivalent to 
\beq
&& \frac{\|\dxp\|^4}{(\pi\delta_L d^{n+1})^2 } \frac{1} {\pi^{(n-k)(k-1)+k}} Vol(\text{Gr}_\C(k-1,n-1))  ... \\
&& ...  \ E^\C_{k-1}(|\det |^{2(n-k+2)}) e_\C (n-k) \big(\int_{\elstar_{|x}\setminus \{0\}}  
\frac{e^{-|a^1_1|^2 } }{\|\lambda\|^{2(k+2)}} d\la   \big) \frac{\omega^n}{n!},
\eeq
where $e_\C(n-k)= \int_{Sym_\C (n-k)}|\det A|^2 d\mu_\C(A)$
and 
$$E^\C_{k-1}(|\det |^{2(n-k+2)})= \int_{M_{k-1}(\C)} |\det A|^{2(n-k+2)}d\mu_\C(A).$$
Now, 
\beq
\frac{\|\dxp\|^4 }{(\pi \delta_L d^{n+1})^2} \int_{\elstar_{|x}\setminus\{0\}}
\frac{e^{-|a^1_1|^2} }{\|\lambda\|^{2k+4}} d\la  &= &
Vol(S^{2k-1}) \frac{\|\dxp\|^4 }{(\pi \delta_L d^{n+1})^2} \int_0^{+\infty} 
\frac{e^{-|a^1_1|^2} }{\|\lambda\|^{5}} d\|\la\|\\  & = &
Vol(S^{2k-1}) \int_0^{+\infty} e^{-r^2} r^3 dr = \frac{1}{2} Vol(S^{2k-1}).
\eeq
Hence, $\E(\nu)$ is equivalent to 
\beq 
 \frac{1} {2\pi^{(n-k)(k-1)+k}} Vol(\text{Gr}_\C(k-1,n-1))  Vol(S^{2k-1}) \ E^\C_{k-1}(|\det |^{2(n-k+2)}) e_\C (n-k) \frac{\omega^n}{n!},
\eeq
where $e_\C (n-k)= (n-k+1)! $ by Proposition 3.8 of \cite{GaWe3},
$Vol(S^{2k-1})=  2\pi^k /(k-1)!,$
$$ E^\C_{k-1}(|\det |^{2(n-k+2)})  =  \frac {\prod_{j=1}^{k-1}\Gamma((n-k+2)+j )}{\prod_{j=1}^{k-1}\Gamma(j)}=
\frac {\prod_{j=n-k+3}^{n+1}\Gamma(j )}{\prod_{j=1}^{k-1}\Gamma(j)}$$
by formula 15.4.12 of \cite{Mehta}
and 
$$Vol(\text{Gr}_\C(k-1,n-1))= \frac{\prod_{j=1}^{k-1} \Gamma(j)}{\prod_{j=n-k+1}^{n-1} \Gamma(j)}\pi^{(k-1)(n-k)}$$
by a computation analogous to the one given
in the real case by Remark \ref{rk311}.
 We conclude that 
 $ \E(\nu)$
 weakly converges to 
${k-1 \choose n-1} \omega^n$,
where the convergence is dominated
by a function in $L^1(X,\frac{\omega^n}{n!})$,
for it has poles of order at most $2(n-k)$ near the critical
points of $p$ and at most $2$ near the base points,
see \cite{GaWe2}. 
\epr
\begin{Corollary}\label{coro cri}
Under the hypotheses of Theorem \ref{lefschetz},
for every generic $\si \in \rhxed$, 
let $|\text{Crit } p_{|\csi} |$ be the number 
of critical points of $p_{|\csi}$. Then,
$$ \frac{1}{d^n } \E(|\text{Crit } p_{|\csi}| )\equid  {k-1 \choose n-1} \int_X c_1(L)^n.$$ 
\end{Corollary}
\bpr 
Corollary \ref{coro cri} follows from Theorem \ref{lefschetz} 
by integration of $1$ over $X$. A direct proof can be given though. 
The modulus of $p$ is a Morse function on $\csi \setminus (\text{Base }(p)\cup F_0\cup F_\infty)$,
where $F_0$ (resp. $F_\infty$) is the fibre of $0$ (resp. of $\infty$) of $p : X \dashrightarrow \cpun$.
Moreover, the index of every critical point of $|p|$ is $n-k$. 
As in the proofs of  Propositions 1 and 2 in \cite{GaWe2}, we deduce that
$\E(|\text{Crit } p_{|\csi} |)$ is equivalent to $|\chi(\csi)|$ as $d$ grows to infinity.
Now, 
$$\chi(\csi)= \int_{\csi} c_{n-k}(\csi) = \int_X c_{n-k}(\csi)\wedge c_k (E\otimes L^d),$$ 
while from the adjunction formula, $c(\csi )\wedge c(E\otimes L^d)_{|\csi}= c(X).$
Moreover, for $0\leq i\leq k$, $c_i (E\otimes L^d) = {k \choose i} d^i c_1(L)^i + o(d^i),$
so that $$c(E\otimes L^d) = (1+dc_1(L))^k + o((1+dc_1(L))^k).$$
From the formula $(1+x)^{-k}= \sum_{j=0}^\infty (-1)^j \frac{(k-1+j)!}{j!(k-1)!} x^j,$
we then deduce that 
$ c_{n-k}(\csi)= (-1)^{n-k} {n-1 \choose k-1} d^{n-k} c_1(L)^{n-k} + o(d^{n-k})$ 
and finally that $$\chi(\csi ) = (-1)^{n-k}  {n-1 \choose k-1} d^{n} \int_X c_1(L)^n+ o(d^n).$$
Hence the result.
\epr

\noindent
\textsc{Damien Gayet\\
UJF-Grenoble 1/CNRS \\
Institut Fourier,}\\
UMR 5582,\\
Grenoble F-38401, France \\
damien.gayet@ujf-grenoble.fr\\

\noindent
\textsc{Jean-Yves Welschinger\\
Universit\'e de Lyon \\
CNRS UMR 5208 \\
Universit\'e Lyon 1 \\
Institut Camille Jordan} \\
43 blvd. du 11 novembre 1918 \\
F-69622 Villeurbanne cedex, France\\
welschinger@math.univ-lyon1.fr

\end{document}